\documentclass[a4paper,12pt,final]{article}
\def\ftype{png}
\pdfoutput=1
\usepackage{geometry}
\usepackage[T1]{fontenc}
\usepackage{ifthen}
\usepackage{amsmath}
\usepackage{amsfonts}
\usepackage{amsthm}
\usepackage{amssymb}
\usepackage[all]{xy}
\usepackage{epsfig}
\usepackage{color}
\usepackage{hyperref}
\hypersetup{colorlinks=true,linkcolor=black,citecolor=black}
\graphicspath{{Graphics/}}
\newtheorem{Lem}{Lemma}[section]
\newtheorem{Pro}[Lem]{Proposition}
\newtheorem*{Con*}{Conjecture}
\newtheorem{Cor}[Lem]{Corollary}
\newtheorem{Thm}[Lem]{Theorem}
\newtheorem*{Thm*}{Theorem}
\theoremstyle{definition}

\newtheorem*{Sch*}{Schedule}
\newtheorem{Exa}[Lem]{Example}
\newtheorem*{Exa*}{Example}
\newtheorem{Def}[Lem]{Definition}
\newtheorem{Rem}[Lem]{Remark}
\newtheorem*{Rem*}{Remark}
\newcommand{\bA}{{\mathbb A}}
\newcommand{\bB}{{\mathbb B}}
\newcommand{\bC}{{\mathbb C}}
\newcommand{\bE}{{\mathbb E}}

\newcommand{\bN}{{\mathbb N}}

\newcommand{\bR}{{\mathbb R}}
\newcommand{\bS}{{\mathbb S}}

\newcommand{\cA}{{\mathcal A}}

\newcommand{\cF}{{\mathcal F}}

\newcommand{\cL}{{\mathcal L}}
\newcommand{\cM}{{\mathcal M}}
\newcommand{\cN}{{\mathcal N}}

\newcommand{\cT}{{\mathcal T}}

\newcommand{\sex}{\operatorname{sup}_{\perp}}
\newcommand{\idty}{{\rm 1\mskip-4mu l}}
\newcommand{\scal}[1]{\langle #1\rangle}
\newcommand{\ri}{\operatorname{ri}}
\newcommand{\rb}{\operatorname{rb}}
\newcommand{\conv}{\operatorname{conv}}
\newcommand{\pos}{\operatorname{pos}}
\newcommand{\nc}{\operatorname{N}}
\newcommand{\tc}{\operatorname{T}}
\newcommand{\tr}{\operatorname{tr}}
\newcommand{\aff}{\operatorname{aff}}
\newcommand{\lin}{\operatorname{lin}}
\begin{document}
\thispagestyle{empty}
\begin{center} 
\textbf{\Large{A Note on Touching Cones and Faces}}\\
\vspace{.3cm}
Stephan Weis\footnote{\texttt{weis@mi.uni-erlangen.de}}\\
Department Mathematik,
Friedrich-Alexander-Universit\"at Erlangen-N\"urnberg,\\
Bismarckstra{\ss}e 1$\frac{\text{1}}{\text{2}}$, D-91054 Erlangen, 
Germany.\\
\vspace{.1cm}
May 3, 2011
\end{center}
\noindent
{\small\textbf{\emph{Abstract --}}
We study touching cones of a (not necessarily closed) convex set in a 
finite-dimensional real Euclidean vector space and we draw relationships 
to other concepts in Convex Geometry. {\it Exposed faces} correspond to 
{\it normal cones} by an antitone lattice isomorphism. {\it Poonems} 
generalize the former to {\it faces} and the latter to
{\it touching cones}, these extensions are non-isomorphic, though. We study
the behavior of these lattices under projections to affine subspaces and 
intersections with affine subspaces. We prove a theorem that characterizes 
exposed faces by assumptions about touching cones. For a convex body $K$ 
the notion of conjugate face adds an isotone lattice isomorphism from the 
exposed faces of the polar body $K^\circ$ to the normal cones of $K$. This 
extends to an isomorphism between faces and touching cones.\\[1mm]
{\em Index Terms\/} -- convex set, exposed face, normal cone, poonem, 
face, touching cone, projection, intersection.\\[1mm]
{\sl AMS Subject Classification:} 52A10, 52A20, 94A17.
%
%
%
%
%
%
%
%
\section{Introduction}
\par
The term of {\it touching cone} has first appeared in 1993 when Schneider 
used it to conjecture\footnote{All of these conjectures are still open.} 
in Section 6.6 of \cite{schneider} equality conditions for the 
Aleksandrov-Fenchel inequality. This inequality, established in 1937, is 
really a system of quadratic inequalities between several
{\it convex bodies}, i.e.\ compact convex subsets of a finite-dimensional 
real Euclidean vector space $(\bE,\langle\cdot,\cdot\rangle)$. A very 
special case is the isoperimetric inequality in dimension two that 
states that the area $A$ and the boundary length $l$ of a two-dimensional 
convex body satisfy $4\pi A\leq l^2$ with equality if and only if the 
convex body is a disk.
\par
Initially we were trying to improve our understanding of projections of 
state spaces. These convex bodies, motivated in Section~\ref{sec:context}, 
are examples where the notion of touching cone is the same as normal cone.
We are not aware of further attention to touching cones in the literature.
So in Section~\ref{subsec:arg} we take the opportunity and collect 
evidence of their significance in Convex Geometry: 
\begin{enumerate}
\item Touching cones arise from normal cones in an analogous way as 
faces arise from exposed faces.
\item
The pair of exposed face and face changes its role with the pair of normal 
cone and touching cone when projection to an affine subspace is replaced 
by intersection with an affine subspace.
\item
If $K$ is a convex body, there is a lattice isomorphism. The faces of the 
polar body correspond to the touching cones of $K$ by taking positive 
hulls.
\item
Touching cones can detect the exposed faces which are intersections of 
coatoms.
\item
Touching cones relate to a special smoothness in dimension two.
\end{enumerate}
\subsection{Preliminaries}
\par
Our analysis uses the frame of Lattice Theory, see e.g.\ Birkhoff 
\cite{birkhoff}, which is well-known in Convex Geometry, see 
e.g.\ Loewy and Tam~\cite{Loewy} and the references therein. A 
mapping $f:X\to Y$ between two partially 
ordered sets (posets) $(X,\leq)$ and $(Y,\leq)$ is {\it isotone} if for 
all $x,y\in X$ such that $x\leq y$ we have $f(x)\leq f(y)$. The mapping 
$f$ is {\it antitone} if for all $x,y\in X$ such that $x\leq y$ we have 
$f(x)\geq f(y)$. A lattice $\cL$ is a partially ordered set $(\cL,\leq) $ 
where the {\it infimum} $x\wedge y$ and {\it supremum} $x\vee y$ of each 
two elements $x,y\in\cL$ exist. All lattices appearing in this article are 
{\it complete}, i.e.\ for an arbitrary subset $S\subset\cL$ the infimum 
$\bigwedge S$ and the supremum $\bigvee S$ exist. The reason is that 
elements $x,y$ in these lattices are convex subsets of $\bE$ where a 
relation $x\leq y$ and $x\neq y$ always implies a dimension step 
$\dim(x)<\dim(y)$ (so $\cL$ has finite length and must be complete). In 
particular $\cL$ has a {\it smallest element} $0$ and a
{\it greatest element} $1$. A {\it coatom} of $\cL$ is an element 
$x\in\cL$ not $1$ such that $y\geq x$ and $y\neq x$ implies $y=1$ for all 
$y\in\cL$.
\begin{figure}[t!]
\centerline{%
\begin{picture}(7.3,3)
\ifthenelse{\equal{\ftype}{eps}}{%
\put(0.2,0.2){\includegraphics[height=2.8cm, bb=0 20 415 400]%
{touching-cones.eps}}}{%
\put(0.2,0.2){\includegraphics[height=2.8cm, bb=0 20 415 400]%
{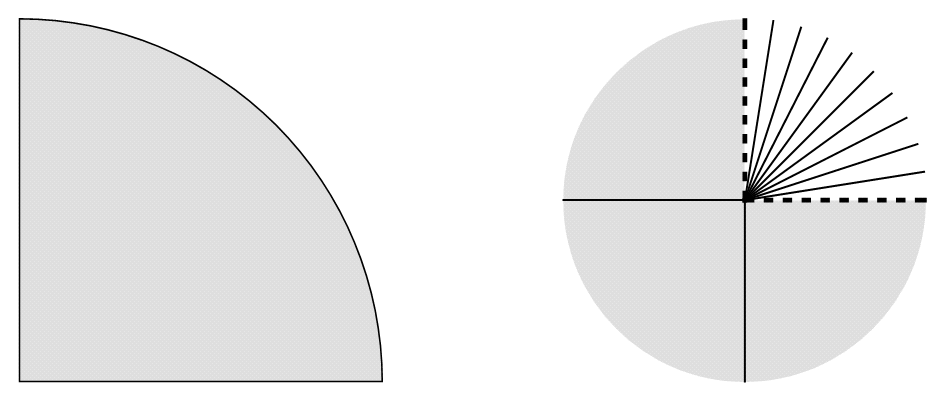}}}
\put(0,0){$a$}
\put(3.2,0){$b$}
\put(0,2.8){$c$}
\put(1,1){$K$}
\end{picture}}
\caption{\label{fig:quarter_disk}The closed quarter disk $K$ (left) with 
its normal cones (right) sketched in the unit disk. Proper normal cones of 
$K$ are: three quadrants at the faces $\{a\},\{b\}$ and $\{c\}$, two rays 
at the faces $[a,b]$ and $[a,c]$ and a family of rays at the one-point 
faces of the arc from $b$ to $c$ other than $\{b\}$ or $\{c\}$. The two 
dashed rays are touching cones but not normal cones of $K$.}
\end{figure}
\par
Given a convex subset $C\subset\bE$ we explain the concepts of normal 
cone, exposed face and face. The {\it normal cone} of $C$ at $x\in C$ is 
the set of vectors $u\in\bE$, that do not make for any $y\in C$ an acute 
angle with the vector from $x$ to $y$. We put
$\nc(C,x):=\{u\in\bE:\scal{u,y-x}\leq 0\text{ for all }y\in C\}$. The
{\it relative interior} $\ri(C)$ of $C$ is the interior of $C$ with 
respect to the affine span $\aff(C)$ of $C$. The {\it relative boundary} 
of $C$ is $\rb(C):=C\setminus\ri(C)$. The {\it normal cone} of any 
non-empty convex subset $F\subset C$ is well-defined (see 
Section~\ref{sec:normal_cones}) as the normal cone of any $x\in\ri(F)$. We 
put $\nc(C,F):=\nc(C,x)$. E.g.\ the normal cone of $C$ is the orthogonal 
complement of the translation vector space $\lin(C)$ of $\aff(C)$ and 
further Examples are shown in Figure~\ref{fig:quarter_disk}. The normal 
cone of the empty set is $\nc(C,\emptyset):=\bE$. This and 
$\lin(C)^\perp$ are the {\it improper normal cones}, all other normal 
cones are {\it proper normal cones} and both together form the {\it normal 
cone lattice} $\cN(C)$. The normal cone lattice is a complete lattice 
ordered by inclusion with the intersection as the infimum (see 
Prop.~\ref{normal_cone_intersection}).%
\begin{figure}[t!]
\centerline{%
\begin{picture}(11,3.3)
\ifthenelse{\equal{\ftype}{eps}}{%
\put(0,0){\includegraphics*[height=3.3cm, bb=0 0 300 200]{stadium.eps}}
\put(6.3,0){\includegraphics[height=3.1cm,bb=0 0 300 205]{mouse.eps}}}{%
\put(0,0){\includegraphics*[height=3.3cm, bb=0 0 380 240]{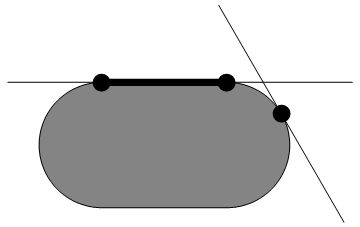}}
\put(6.3,0){\includegraphics[height=3.1cm,bb=0 0 370 250]{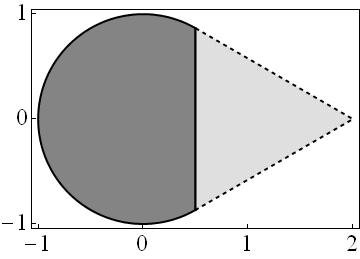}}}
\put(0.2,2.2){$H_1$}
\put(2.1,2.2){$F_1$}
\put(4.8,0.2){$H_2$}
\put(4.1,1.4){$F_2$}
\put(7.5,1.5){\textcolor{white}{$K$}}
\end{picture}}
\caption{\label{fig:exposed_faces_normal_cones}%
The stadium (left) consists of a square with two half-disks attached on 
opposite sides. The supporting hyperplane $H_i$ defines the exposed 
face $F_i$ for $i=1,2$. The two extreme points of $F_1$ are non-exposed 
faces. The truncated disk $K$ (right) is the closed unit ball in $\bR^2$ 
with the segment $x>\frac{1}{2}$ missing. The polar body 
$K^\circ$ of $K$ is the union of $K$ with the bright closed triangle.}
\end{figure}
\par 
A {\it supporting hyperplane} of $C$ is any affine hyperplane $H$ in 
$\bE$, such that $C\setminus H$ is convex and $C\cap H$ is non-empty. An 
{\it exposed face} of $C$ is the intersection of $C$ with a supporting 
hyperplane. An example is shown in 
Figure~\ref{fig:exposed_faces_normal_cones}, left. In addition 
$\emptyset$ and $C$ are exposed faces called {\it improper exposed faces}. 
All other exposed faces are {\it proper exposed faces}. The set of exposed 
faces is the {\it exposed face lattice} $\cF_\perp(C)$. This is a complete 
lattice ordered by inclusion and with the intersection as the infimum (see 
Prop.~\ref{pro:ex_lattice}). If $C$ has at least two points, then we have
an antitone lattice isomorphism (see Prop.~\ref{lattice_antitone_iso})
\begin{equation}
\label{eq:anti_iso}
\nc(C):\quad\cF_\perp(C)\to\cN(C),\quad
F\mapsto\nc(C,F).
\end{equation}
Two examples of this isomorphism are sketched in Figure~\ref{fig:lattices}
in columns two and three. The isomorphism does not require that $C$ is 
closed or bounded. We can write the isomorphism
(\ref{eq:anti_iso}) and its inverse in the form (\ref{eq:anti_isos}), i.e.\
for proper exposed faces $F$ and proper normal cones $N$ of $C$ we have
\[\textstyle
\begin{array}{rcl}
F & \mapsto & \bigcap_{x\in F}\nc(C,x)
\quad = \quad \nc(C,y) \qquad\text{for any }y\in\ri(F),\\
N & \mapsto & \bigcap_{u\in N\setminus\{0\}}F_\perp(C,u)
\quad = \quad F_\perp(C,v) \qquad\text{for any }v\in\ri(N)\setminus\{0\}.
\end{array}\]
\par
The {\it closed segment} between $x,y\in\bE$ is
$[x,y]:=\{(1-\lambda)x+\lambda y\mid\lambda\in[0,1]\}$, the 
{\it open segment} between $x,y\in\bE$ is
$]x,y[\,:=\{(1-\lambda)x+\lambda y\mid\lambda\in(0,1)\}$. A {\it face} of 
$C$ is a convex subset $F$ of $C$, s.t.\ whenever for 
$x,y\in C$ the open segment $]x,y[$ intersects $F$, then the closed segment
$[x,y]$ is included in $F$. An {\it extreme point} is the element 
of a zero-dimensional face. The faces $\emptyset$ and $C$ are
{\it improper faces}, all other faces are {\it proper faces}. The set of 
all faces of $C$ is the {\it face lattice} of $C$ denoted by $\cF(C)$. It 
is easy to show that arbitrary intersections of faces are faces, so 
$\cF(C)$ is a complete lattice ordered by inclusion and with the 
intersection as the infimum. It is easy to show
$\cF(C)\supset\cF_\perp(C)$. A face which is not an exposed face will be 
called a {\it non-exposed} face, see e.g.\ 
Figure~\ref{fig:exposed_faces_normal_cones}, left. 
\subsection{Observations about touching cones}
\label{subsec:arg}
\par
We introduce touching cones according to our results in 
Theorem~\ref{thm:rtc_cover}. A {\it touching cone} of $C$ is any non-empty 
face of a normal cone of $C$. An example is shown in 
Figure~\ref{fig:quarter_disk}. The improper normal cones $\lin(C)^\perp$ 
and $\bE$ are touching cones called {\it improper touching cones}, all 
other touching cones are {\it proper touching cones}. These together form 
the {\it touching cone lattice} denoted by $\cT(C)$. This is a complete 
lattice ordered by inclusion and with the intersection as the infimum. 
One has $\cT(C)\supset\cN(C)$.
\subsubsection{Analogy in creation touching cones and faces}
\par
There is an analogy between touching cone and face if we use the 
concept of {\it poonem} that Gr\"unbaum \cite{gruenbaum} applies 
for a closed convex subset of $\bE$. In finite dimension {\it poonem} is 
equivalent to {\it face}. We define a {\it poonem} of a convex subset 
$C\subset\bE$ as a subset $P$ of $C$ s.t.\ there exist subsets 
$F_0,F_1,\ldots,F_k$ of $C$ with $F_0=P$, $F_k=C$ and $F_{i-1}$ is an 
exposed 
face of $F_i$ for $i=1,\ldots,k$. Every poonem is a face because a face of 
a face of $C$ is a face of $C$. The converse is also true: given a proper 
face $F$ of $C$, the smallest exposed face $\sex(F)$ containing $F$ is a 
proper exposed face of $C$ by Lemma~\ref{lem:sup_exp}, so 
$\dim(\sex(F))<\dim(C)$. By induction $F$ is a poonem of $C$. We have 
unified extensions
\[\begin{array}{rcl}
\cF_\perp(C) & \subset & \cF(C)=
\{\text{ poonems of elements in }\cF_\perp(C)\},\\
\cN(C) & \subset & \cT(C)=
\{\text{ non-empty poonems of elements in }\cN(C)\}.
\end{array}
\]
As $\cF(C)$ is the set of poonems of $C$, a more 
systematic definition would consider poonems of proper elements or of 
coatoms of $\cF_\perp(C)$ and of $\cN(C)$. In any case we can see that the 
concepts of {\it exposed face}, {\it normal cone} and {\it poonem} suffice 
to define {\it face} and {\it touching cone} in a unified way.
\subsubsection{Compatibility with projection and intersection}
\par
We introduce Schneider's (equivalent) definition of touching cone:
If $v\in\bE$ is non-zero and the exposed face $F:=F_\perp(C,v)$ is 
non-empty, then the face $T(C,v)$ of the normal cone $\nc(C,F)$ that 
contains $v$ in its relative interior, is called a touching cone; 
$\lin(C)^\perp$ and $\bE$ are touching cones by definition.
\par
Let $\bA\subset\bE$ be an affine subspace, by $\pi_\bA(C)$ we denote the 
orthogonal projection of $C$ to $\bA$. If $v\in\lin(\bA)$ and $T(C,v)$ is 
a normal cone of $C$, then $T(\pi_\bA(C),v)$ is a normal cone of 
$\pi_\bA(C)$. This is proved in Section \ref{sec:sharp} by a new 
characterization of normal cones. Exposed faces of $C$ however may 
project to non-exposed faces of $\pi_\bA(C)$.
\par
Dually, exposed faces are preserved under intersection of $C$ with $\bA$. 
But for some $v\in\lin(\bA)$ the cone $T(C,v)$ may be a normal cone of 
$C$ while $T(C\cap\bA,v)$ is not  a normal cone of $C\cap\bA$. 
Example~\ref{ex:ip} discusses these aspects.
\subsubsection{A lattice isomorphism for convex bodies}
\par
We consider a convex body $K\subset\bE$ with at least two points 
and with the origin in the interior, $0\in{\rm int}(K)$. The
{\it polar body}
\[
K^\circ:=\{u\in\bE\mid\langle u,x\rangle\leq 1\text{ for all }x\in K\}
\]
is a convex body with $0\in{\rm int}(K^\circ)$, an example is shown in 
Figure~\ref{fig:exposed_faces_normal_cones}, right. Given a subset 
$S\subset\bE$, the {\it positive hull} $\pos(S)$ of $S$ is the set of all 
finite {\it positive combinations} of elements of $S$, i.e.\ an element 
$x\in\bE$ belongs to $\pos(S)$ if and only if there is $k\in\bN$, 
$\lambda_i\in\bR$ with $\lambda_i\geq 0$ and $s_i\in S$ for $i=1,\ldots,k$ 
such that $x=\sum_{i=1}^k\lambda_is_i$ (we have $0\in\pos(S)$). In 
Section~\ref{sec:polarity} we establish isotone lattice isomorphisms 
\begin{equation}
\label{eq:isot_iso}
\begin{array}{rcl}
\cF_\perp(K^\circ)\to\cN(K), && F\mapsto \pos(F),\\
\cF(K^\circ)\to\cT(K), && F\mapsto \pos(F).
\end{array}
\end{equation}
The inverse isomorphism is given for a proper touching cone $T\in\cT(K)$ by
$T\mapsto\rb(K^\circ)\cap T$. We think that (\ref{eq:isot_iso}) underlines 
(in the case of convex bodies) that the notion of touching cone is as 
fundamental as face. An example of the lattice isomorphisms is shown in 
Figure~\ref{fig:lattices}.%
\begin{figure}[t!]
\centerline{%
\begin{picture}(15,5.8)
%
\ifthenelse{\equal{\ftype}{eps}}{%
\put(0,3.4){\includegraphics[width=3cm,bb=0 0 497 400]{C1.eps}}
\put(4,3.4){\includegraphics[width=3cm,bb=0 0 497 400]{C11.eps}}
\put(8,3.4){\includegraphics[width=3cm,bb=0 0 497 400]{C12.eps}}
\put(12,3.4){\includegraphics[width=3cm,bb=0 0 497 400]{C13.eps}}
\put(0,0){\includegraphics[width=3cm,bb=0 0 497 400]{C2.eps}}
\put(4,0){\includegraphics[width=3cm,bb=0 0 497 400]{C21.eps}}
\put(8,0){\includegraphics[width=3cm,bb=0 0 497 400]{C22.eps}}
\put(12,0){\includegraphics[width=3cm,bb=0 0 497 400]{C23.eps}}}{%
\put(0,3.4){\includegraphics[width=3cm,bb=0 0 497 400]{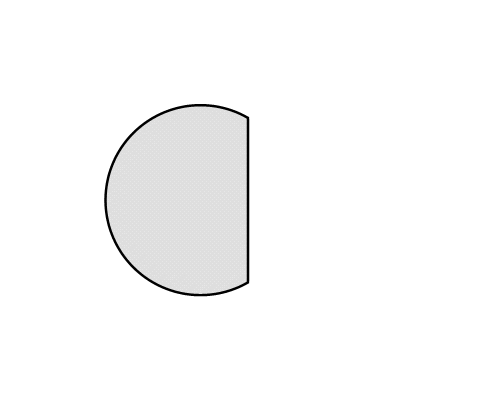}}
\put(4,3.4){\includegraphics[width=3cm,bb=0 0 497 400]{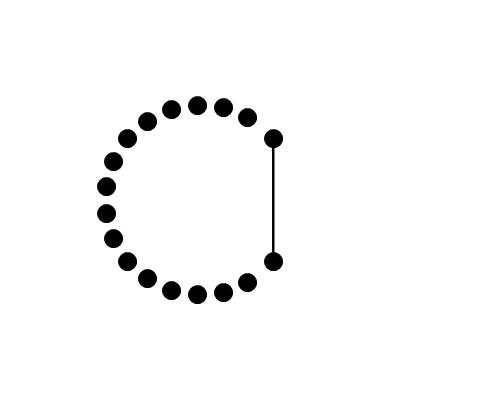}}
\put(8,3.4){\includegraphics[width=3cm,bb=0 0 497 400]{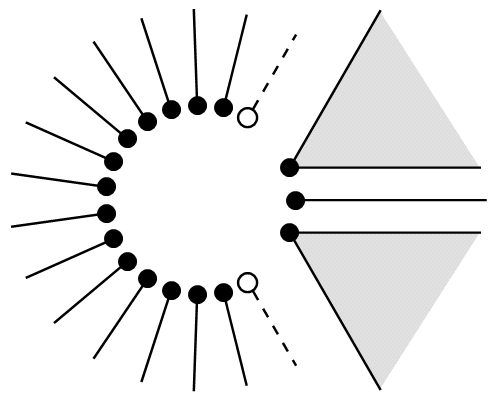}}
\put(12,3.4){\includegraphics[width=3cm,bb=0 0 497 400]{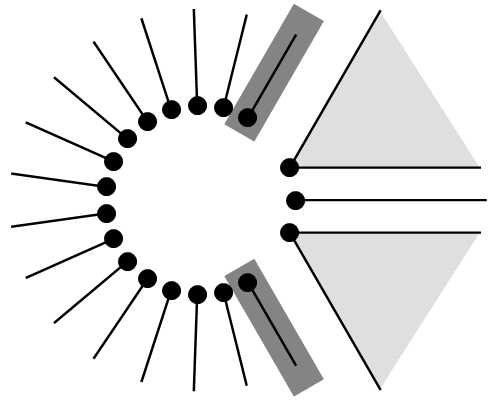}}
\put(0,0){\includegraphics[width=3cm,bb=0 0 497 400]{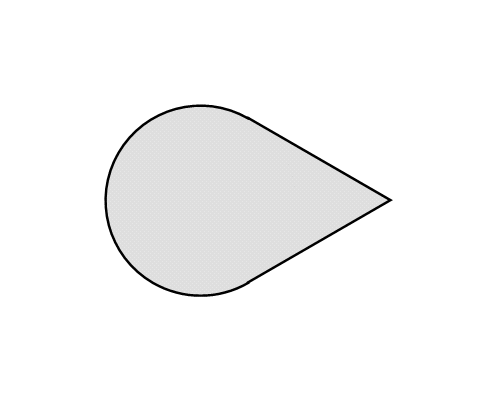}}
\put(4,0){\includegraphics[width=3cm,bb=0 0 497 400]{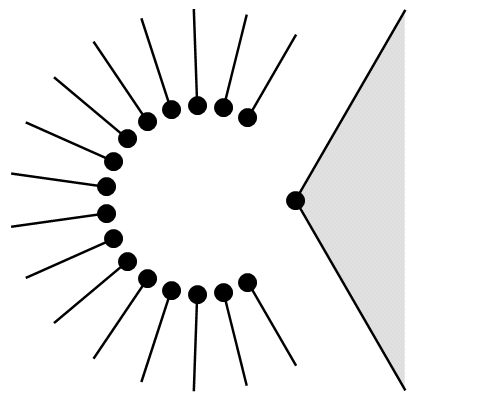}}
\put(8,0){\includegraphics[width=3cm,bb=0 0 497 400]{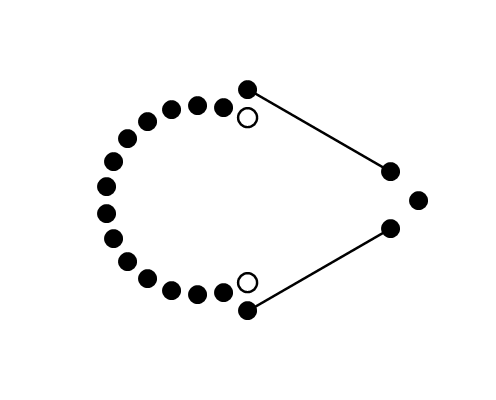}}
\put(12,0){\includegraphics[width=3cm,bb=0 0 497 400]{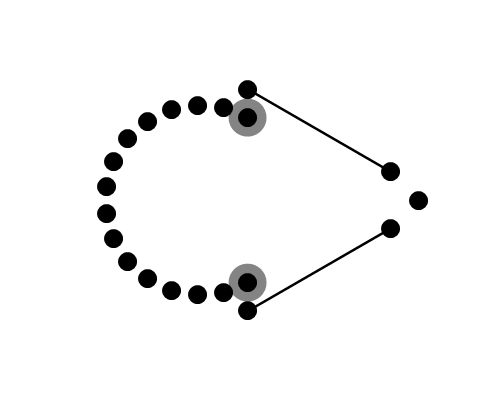}}}
\put(1.0,4.5){$K$}
\put(1.7,4.5){$=K^{\circ\circ}$}
\put(5.0,4.5){$\cF_{\perp}$}
\put(9.0,4.5){$\cN$}
\put(13.0,4.5){$\cT$}
\put(7.0,4.5){$\stackrel{\nc(K)}{\longrightarrow}$}
\put(11.3,4.5){$\subset$}
\put(1.0,1.1){$K^{\circ}$}
\put(5.0,1.1){$\cN$}
\put(9.0,1.1){$\cF_{\perp}$}
\put(13.0,1.1){$\cF$}
\put(7.0,1.1){$\stackrel{\nc(K^{\circ})}{\longleftarrow}$}
\put(11.3,1.1){$\subset$}
\put(1.4,2.8){$\Big\updownarrow\text{polar body}$}
\put(5.4,2.8){$\Big\downarrow\pos$}
\put(7.4,2.7){$\!\!\!\nwarrow\!\!\searrow^{F\mapsto\widehat{F}}$}
\put(9.4,2.8){$\Big\uparrow\pos$}
\put(13.4,2.8){$\Big\uparrow\pos$}
\end{picture}}
\caption{\label{fig:lattices}A finite sketch of proper lattice 
elements, empty circles denote deleted points, dashed lines denote deleted 
lines. Lattices belong to the convex body to their left, we have
$\cF(K)=\cF_\perp(K)$ and $\cT(K^\circ)=\cN(K^\circ)$. In both 
rows there is an antitone isomorphism between exposed 
faces and normal cones (between columns two and three). The positive hull 
operator $\pos$ defines three isotone isomorphisms between rows one and 
two. Touching cones that are not normal cones and non-exposed faces are 
highlighted by a dark background (right column). The antitone isomorphism 
of the conjugate face is $F\mapsto\widehat{F}$.}
\end{figure}
\par
Following Remark~\ref{rem:tc_cover} for a convex body $K$ we have the 
partition of $\bE$ into the relative interiors of touching cones 
$\neq\bE$. Denoting $T(K,u)$ the touching cone with the vector 
$u\in\bE\setminus\{0\}$ in its relative interior, we have
the partition
\[\textstyle
\bE\setminus\{0\}=\bigcup\limits^{\bullet}{}_{u\in\bE\setminus\{0\}}
\ri(T(K,u)).
\]
This is reminiscent of the partition of the {\it metric projection} 
(see e.g.\ Schneider \cite{schneider})
\[\textstyle
\bE=\bigcup\limits^{\bullet}{}_{x\in K}(x+\nc(K,x)).
\]
The partition of $\bE\setminus\{0\}$ reminds us also of the partition of 
$K^\circ$ into the relative interiors of its faces (\ref{eq:stratum}). We 
have the following analogy:
\[\begin{array}{l|l}
\text{Partition of }\rb(K^\circ)\text{ in relative interiors}
& \text{Partition of }\bE\setminus\{0\}\text{ in relative interiors}\\
\text{of proper faces.} & \text{of proper touching cones of }K.
\end{array}\]
\subsubsection{Coatoms of the face lattice}
\par
We explain for a general convex subset $C\subset\bE$ that touching 
cones can characterize exposed faces in terms of coatoms in 
$\cF_\perp(C)$. We recall that a coatom $F$ of $\cF_\perp(C)$ does not 
need to satisfy the dimension equation $\dim(F)+1=\dim(C)$, see e.g.\ 
$F_2$ in Figure~\ref{fig:exposed_faces_normal_cones}, left. Since 
intersections of exposed faces are exposed, any intersection of coatoms in 
$\cF_\perp(C)$ is an exposed face. A sufficient condition for the converse 
is proved in Thm.~\ref{thm:local}:
\begin{Thm*}
Let $F$ be a proper exposed face of $C$ where every touching cone included 
in the normal cone $\nc(C,F)$ is a normal cone. Then $F$ is an 
intersection of coatoms of $\cF_\perp(C)$.
\end{Thm*}
\par
Figure~\ref{fig:lens} shows that there is no converse to the theorem. 
Examples are discussed after the remark below. A main argument to the 
theorem is {\it Minkowski's theorem} (a convex body is the convex hull of 
its extreme points) applied to a section of a normal cone. Another 
argument is the isomorphism (\ref{eq:anti_iso}). If we consider convex 
bodies, then the isomorphism (\ref{eq:isot_iso}) turns the theorem into an 
equivalent form, which more obviously follows from Minkowski's theorem 
(see Section~\ref{sec:polarity}).%
\begin{figure}[t!]
\centerline{%
\begin{picture}(6.2,1.5)
\ifthenelse{\equal{\ftype}{eps}}{%
\put(0,0){\includegraphics[height=1.5cm, bb=5 0 300 70]{lens.eps}}
}{%
\put(0,0){\includegraphics[height=1.5cm, bb=5 0 300 85]{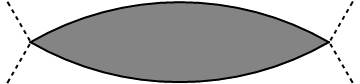}}
}
\end{picture}}
\caption{\label{fig:lens}
This intersection of two closed disks has touching cones which are 
not normal cones (dashed rays). Still, all proper faces are coatoms.}
\end{figure}
\begin{Rem}[Exposed faces in dimension two]
\label{ref:exposed_2D}
In dimension $\dim(C)=2$ every non-exposed face of $C$ is the endpoint of 
a unique one-dimensional face of $C$. 
\par
We prove this claim. All one-dimensional faces of $C$ are coatoms of 
$\cF_\perp(C)$ (as $\sex(F)$ is proper for a proper face $F$). One 
dimension below, a
point $x$ of $C$ may belong to $i=0,1,2$ one-dimensional faces of $C$ and 
exactly for $i=0,2$ the set $\{x\}$ is an intersection of coatoms of $C$.  
So a proper exposed face $F$ of $C$ is not the intersection of coatoms of 
$\cF_\perp(C)$ if and only if $F=\{x\}$ where $x$ is the endpoint of 
a unique one-dimensional face of $C$. 
\par
If in addition the assumptions of the above theorem hold for $C$, then 
non-exposed faces $F$ are characterized by the conditions $F=\{x\}$ where 
$x$ is the endpoint of a unique one-dimensional face of $C$.
\end{Rem}
\begin{figure}[t!]
\centerline{%
\ifthenelse{\equal{\ftype}{eps}}{%
\includegraphics[height=3.2cm,bb=0 0 300 272]{B064sharp_non_closed.eps}
\hspace{2cm}
\includegraphics[height=3.2cm,bb=0 0 300 272]{B065non_sharp.eps}}{%
\includegraphics[height=3.2cm,bb=0 0 360 320]{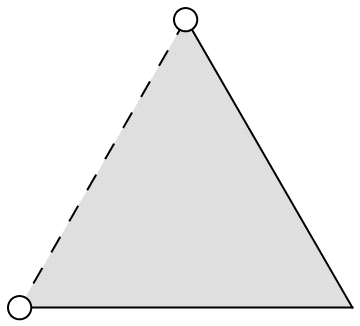}
\hspace{2cm}
\includegraphics[height=3.2cm,bb=0 0 360 320]{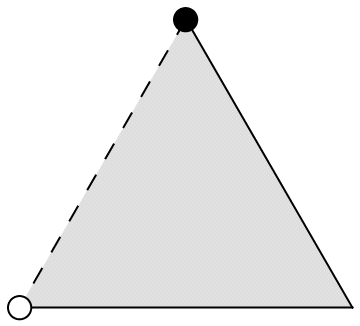}}}
\caption{\label{fig:sharp_non-closed}Empty circles denote deleted 
points, dashed lines denote deleted lines. The left triangle has three 
proper touching cones, all of which are normal cones. Accordingly every 
proper exposed face of the triangle is an intersection of coatoms. If 
the top vertex is added (right) one normal cone is added but two touching 
cones are added. The top vertex is not an intersection of coatoms.}
\end{figure}
\par
An example with $\cN(C)=\cT(C)$ is the polar body $K^\circ$ (mouse shape)
of the truncated disk. Further examples of $\cN(C)=\cT(C)$ are the state 
space discussed in Example~\ref{ex:ip}.
Examples of $\cN(C)\subsetneq\cT(C)$ that do not have the characterization 
of Remark~\ref{ref:exposed_2D} are the quarter disk in 
Figure~\ref{fig:quarter_disk} and the truncated disk $K$ in 
Figure~\ref{fig:lattices}. Two convex set, which are not closed, are 
discussed in Figure~\ref{fig:sharp_non-closed}. 
\subsubsection{Smoothness in dimension two}
\par
There is a special smoothness issue in dimension two. This holds for
a general convex subset $C\subset\bR^2$ if $\cN(C)=\cT(C)$, examples are 
listed in the previous paragraph. It would be 
interesting to see how smoothness generalizes into higher dimensions 
(where however coatoms of $\cF_\perp(C)$ can have small dimension). A 
boundary point $x$ of $C$ is {\it singular}, if $C$ has two linearly 
independent normal vectors at $x$. 
\par
The smoothness property, given $\dim(C)=2$ and $\cT(C)=\cN(C)$, is that 
every singular point $x\in C$ is the intersection to two distinct boundary 
segments of $C$: If $x\in C$ is singular then the normal cone of $C$ at 
$x$ has two distinct boundary rays $t_1,t_2$, which are touching cones of 
$C$ by definition. By assumption $t_1$ is a normal cone of $C$, so it is 
the normal cone at a boundary point $y_1\neq x$ of $C$. It follows that 
the segment $[x,y_1]$ is a boundary segment of $C$. The same arguments 
applied to $t_2$ show $\{x\}=[x,y_1]\cap[x,y_2]$ (If the intersection was 
a segment, then $\dim(C)\leq 1$ by (\ref{eq:nc_for_begn}) (iv)).
\subsection{Projections of state spaces}
\label{sec:context}
\par
Our motivation to study touching cones lies in Information Theory, see 
Amari and Nagaoka \cite{amari_nagaoka}. Analysis takes place in the convex 
body of {\it state space} $\bS(n)$. This is a convex body in the algebra 
${\rm Mat}(\bC,n)$ of complex $n\times n$-matrices. In fact $\bS(n)$ 
consists of all positive semi-definite matrices (i.e.\ being self-adjoint 
and without negative eigenvalues) that have trace one. We have
$\cT(\bS(n))=\cN(\bS(n))$ and $\cF(\bS(n))=\cF_\perp(\bS(n))$, see 
Example~\ref{ex:state_spaces}. In Example~\ref{ex:ip} we discuss
orthogonal projections $P$ of $\bS(n)$ to vector spaces $L$, they too 
satisfy $\cT(P)=\cN(P)$. These projections are connected to information 
manifolds called {\it exponential families}, see e.g.\ Knauf and Weis 
\cite{knauf_weis}.
\par
We ask if a finite-dimensional convex set $C$ is {\it stable}, which means 
that for any $0\leq d\leq\dim(C)$ the union of faces $F$ of $C$ with 
$\dim(F)\leq d$ is a closed set (see Papadopoulou \cite{papadopoulou}). It 
is well-known that $\bS(n)$ is stable. Is $P$ also stable? This would have 
consequences for the topology of exponential families.
\par
Another question is about non-exposed faces of $P$ and their behavior if 
$L$ varies in a Grassmannian manifold of subspaces. This question may be 
related to continuity properties of information measures, see 
\cite{knauf_weis}. It is likely to be accessible by Convex 
Algebraic Geometry (as studied by Henrion, Rostalski, Sturmfels 
and others) because $P$ is polar to an affine section of $\bS(n)$, see 
\cite{henrion,rostalski,Weis}. On the other hand, the faces of $P$ 
correspond to the touching cones of the affine section, which is an affine 
algebraic set.
%
%
%
%
\section{Posets and lattices}
\label{sec:posets}
\par
We introduce lattices and cite two fundamental assertions about lattices. 
\begin{Def}
\label{lattice_def}
A {\it partially ordered set} or {\it poset} $(X,\leq)$ is a set $X$ with 
a binary relation $\leq$, such that for all $x,y,z\in X$ we have
$x\leq x$ (reflexive), $x\leq y$ and $y\leq x$ implies $x=y$ 
(antisymmetric) and $x\leq y$ and $y\leq z$ implies $x\leq z$ 
(transitive); $y\geq x$ is used instead of $x\leq y$.
\par
A mapping $f:X\to Y$ between two posets $(X,\leq)$ and $(Y,\leq)$ is 
{\it isotone}, if $x_1\leq x_2$ implies $f(x_1)\leq f(x_2)$ for any 
$x_1,x_2\in X$. The mapping $f$ is {\it antitone} if $x_1\leq x_2$ implies
$f(x_2)\leq f(x_1)$. 
\par
In a poset $(X,\leq)$, a {\it lower bound} of a subset $S\subset X$ is 
an element $x\in X$ such that $x\leq s$ for all $s\in S$. An
{\it infimum} of $S$ is a lower bound $x$ of $S$ such that $y\leq x$ for 
every lower bound $y$ of $S$. Dually, an {\it upper bound} of a subset 
$S\subset X$ is an element $x\in X$ such that $s\leq x$ for all $s\in S$. 
A {\it supremum} of $S$ is an upper bound $x$ of $S$ such that $x\leq y$ 
for every upper bound $y$ of $S$. We may write 
$S=\{s_{\alpha}\}_{\alpha\in I}$ for an index set $I$. In case of 
existence, the infimum of $S$ is unique and is denoted by $\bigwedge S$ 
or by $\bigwedge_{\alpha\in I}s_{\alpha}$, likewise the supremum of $S$ 
is denoted by $\bigvee S$ or by $\bigvee_{\alpha\in I}s_{\alpha}$ in case 
of existence.
\par
If $(X,\leq)$ has a smallest element $0$, then an element $x\in X$ not $0$ 
is an {\it atom} of $X$ if for all $y\leq x$ in $X$ with $y\neq x$ we have 
$y=0$. If $(X,\leq)$ has a greatest element $1$, then an element $x\in X$ 
not $1$ is a {\it coatom} of $X$ if for all $y\geq x$ in $X$ with
$y\neq x$ we have $y=1$. 
\par
A {\it lattice} $(\cL,\leq,\wedge,\vee)$ is a poset $(\cL,\leq)$, such 
that for any two elements  $x,y\in\cL$ the infimum
$x\wedge y:=\bigwedge\{x,y\}$ and the supremum $x\vee y:=\bigvee\{x,y\}$ 
exist. A lattice $(\cL,\leq,\wedge,\vee)$ is {\it complete} if 
every subset $X$ of $\cL$ has an infimum and a supremum. We denote a 
complete lattice by $(\cL,\leq,\wedge,\vee,0,1)$ with $0$ the smallest and 
$1$ the greatest element of $\cL$. A lattice $(\cL,\leq,\wedge,\vee)$ is 
{\it modular} if for all elements $x,y,z\in\cL$ the {\it modular law} is 
true:
\begin{equation}
\label{eq:mod_l}
x\leq z\quad\text{implies}\quad
x\vee(y\wedge z)=(x\vee y)\wedge z.
\end{equation}
The partial ordering of $\cL$ restricts to subsets. We call 
$X\subset\cL$ a {\it sublattice} of $\cL$ if for all $x,y\in X$ the 
infimum $x\wedge y$ and the supremum $x\vee y$ (calculated in $\cL$) 
belong to $X$.
\end{Def}
\begin{Rem}
\label{rem:isotone_isomorphism}
Birkhoff has proved in \cite{birkhoff}, Lemma 1 on page 24, that
an isotone bijection between two lattices with isotone inverse is a 
lattice isomorphism.
\end{Rem}
\begin{Def}
\label{def:closure_property}
A property of subsets of a set $M$ is a {\it closure property}
when (i) $M$ has the property, and (ii) 
any intersection of subsets having the given property itself has this 
property.
\end{Def}
\begin{Rem}
\label{rem:closure_complete}
Birkhoff has proved in \cite{birkhoff}, Corollary on page 7, that those 
subsets $\cM$ of any set $M$ which have a given closure property form a 
complete lattice. The ordering on $\cM$ is given by inclusion. The 
infimum of $\{M_{\alpha}\}_{\alpha\in I}\subset\cM$ is the intersection
$\bigwedge_{\alpha\in I} M_{\alpha}=\bigcap_{\alpha\in I} M_{\alpha}$
and the supremum is
$\bigvee_{\alpha\in I}M_{\alpha}=\bigcap\{\widetilde{M}\in\cM\mid
\forall \alpha\in I:M_{\alpha}\subset\widetilde{M}\}$. 
\end{Rem}
%
%
%
%
%
%
%
%
\section{Faces and exposed faces}
\label{sec:faces}
\par
We introduce faces and exposed faces of a convex set and their lattice 
structure. Klingenberg \cite{klingenberg} may be consulted for the 
background in affine geometry. Let $(\bE,\langle\cdot,\cdot\rangle)$ be a
finite-dimensional real Euclidean vector space. We recommend a monograph 
such as Rockafellar or Schneider \cite{rock,schneider} for an introduction 
to convex sets.
\begin{Def}[Convexity]
The {\it convex hull} $\conv(C)$ of a subset $C\subset\bE$ consists of all 
{\it convex combinations} of elements of $C$, i.e.\ $x\in\conv(C)$ if and 
only if there is $k\in\bN$ and for $i=1,\ldots,k$ there are 
$\lambda_i\in\bR$ with $\lambda_i\geq 0$ and $\sum_{j=1}^k\lambda_j=1$ and 
there are $x_i\in C$ such that $x=\sum_{j=1}^k\lambda_jx_j$. We understand 
$\conv(\emptyset)=\emptyset$. The subset $C\subset\bE$ is {\it convex}, if 
$x,y\in C$ implies $[x,y]\subset C$, which is the same as $C=\conv(C)$.
A {\it convex body} is a closed and bounded convex set. If we drop the 
condition of $\sum_{i=1}^n\lambda_i=1$ then we speak of a
{\it positive combination} and we denote the set of positive combinations 
of $C$ by $\pos(C)$ (and $\pos(\emptyset)=\{0\}$). A {\it convex cone} is 
a non-emepty convex subset $C$ of $\bE$ where $x\in C$ and $\lambda\geq 0$ 
imply $\lambda x\in C$, which is the same as $C=\pos(C)$. 
\end{Def}
\par
According to Rockafellar \cite{rock} \textsection 2 the convex hull 
of $C$ is the smallest convex subset of $\bE$ containing $C$.  It is a 
closure property that a subset $C\subset\bE$ is 
convex, i.e.\ $\bE$ is convex and arbitrary intersections of convex 
subsets are convex. Hence, Remark~\ref{rem:closure_complete} ensures that 
the convex subsets of $\bE$ are the elements of a complete lattice ordered 
by inclusion and $\conv(C)$ is the intersection of all convex subsets of 
$\bE$ that include $C$. Closure properties are important also for face 
lattices.
\begin{Def}[Face lattice]
\label{def:face_lattice}
If $C\subset\bE$ is a convex subset, then a convex subset $F\subset C$ is 
a {\it face} of $C$ if for all $x,y\in C$ the non-empty intersection 
$]x,y[\,\cap F$ implies  $[x,y]\subset F$. The empty set $\emptyset$ and 
$C$ are {\it improper} faces, all other faces of $C$ are {\it proper}. A 
face of the form $\{x\}$ for $x\in C$ is called an {\it  extreme point} 
of $C$. The set of faces of $C$ will be denoted by $\cF(C)$ and will be 
called the {\it face lattice} of $C$.
\end{Def}
\par
If $C\subset\bE$ is a convex subset then the intersection of any family of 
faces of $C$ is a face of $C$. In other words, the property face is a 
closure property. Thus, by Remark~\ref{rem:closure_complete} the face 
lattice
\begin{equation}
\label{face_lattice}
(\cF(C),\subset,\cap,\vee,\emptyset,C)
\end{equation}
is a complete lattice ordered by inclusion and the infimum is the 
intersection. The smallest element of $\cF(C)$ is $\emptyset$, the 
greatest is $C$. We cite Schneider \cite{schneider}, Chap.\ 1, for
two fundamental theorems. {\it Carath\'eodory's theorem} says if 
$C\subset\bE$ and $x\in\conv(C)$, then $x$ is a convex combination of 
affinely independent points of $C$. {\it Minkowski's theorem} says that 
every convex body is the convex hull of its extreme points. 
\begin{Def}[Relative interior]
If $C\subset\bE$ then the {\it affine hull} of $C$, denoted by $\aff(C)$ 
is the smallest affine subspace of $\bE$ that contains $C$. The interior 
of $C$ with respect to the relative topology of $\aff(C)$ is the
{\it relative interior} $\ri(C)$ of $C$. The complement 
$\rb(C):=C\setminus\ri(C)$ is the {\it relative boundary} of $C$. If 
$C\subset\bE$ is convex and non-empty then the {\it vector space} of $C$ 
is defined as the translation vector space of $\aff(C)$,
\begin{equation}
\label{def:trans_space}
\lin(C):=\{x-y\mid x,y\in\aff(C)\}.
\end{equation}
We define the {\it dimension} $\dim(C):=\dim(\lin(C))$ and 
$\dim(\emptyset)=-1$.
\end{Def}
\par
Let $C,D\subset\bE$ be convex subsets. Rockafellar proves in \cite{rock}, 
Coro.~6.6.2, the sum formula for the 
relative interior
\begin{equation}
\label{eq:risum}
\ri(C)+\ri(D)=\ri(C+D).
\end{equation}
In Thm.\ 6.5 he proves for the case $\ri(C)\cap\ri(D)\neq\emptyset$
\begin{equation}
\label{eq:riint}
\ri(C)\cap\ri(D)=\ri(C\cap D).
\end{equation}
If $\bA$ is an affine space and $\alpha:\bE\to\bA$ is an affine mapping, 
then by Thm.\ 6.6 in 
\cite{rock}
\begin{equation}
\label{eq:cori}
\alpha(\ri(C))=\ri(\alpha(C))
\end{equation}
holds. If $F$ is a face of $C$ and if $D$ is a subset of $C$, then
by Thm.\ 18.1 in \cite{rock} we have
\begin{equation}
\label{eq:ri_inc}
\ri(D)\cap F\neq\emptyset \implies  D\subset F.
\end{equation}
By  Thm.\ 18.2 in \cite{rock} $C$ admits a partition by relative interiors 
of its faces
\begin{equation}\textstyle
\label{eq:stratum}
C=\bigcup\limits^{\bullet}{}_{F\in\cF(C)}\ri(F).
\end{equation}
In particular, every proper face of $C$ is included in the relative 
boundary $\rb(C)$ and its dimension is strictly smaller than the dimension 
of $C$. We need the following.
\begin{Lem}
\label{lem:cone_faces}
If $H\subset\bE$ is an affine hyperplane with $0\not\in H$ and $C\subset H$
is a convex subset, then 
$\pos:\cF(C)\to\cF(\pos(C))\setminus\{\emptyset\}$ is a bijection with
inverse $F\mapsto C\cap F$.
\end{Lem}
{\em Proof:\/}
If $F$ is a face of $\pos(C)$, then $F$ is a convex cone. So, if 
$F\neq\emptyset$, then $F=\pos(F\cap C)$. Moreover, since 
$C\subset\pos(C)$ the set $F\cap C$ is a face of $C$. This gives an 
injective mapping
\[
\cF(\pos(C))\setminus\{\emptyset\}\to\cF(C),\quad
F\mapsto F\cap C.
\]
By (\ref{eq:stratum}) the relative interiors of faces $F$ of $\pos(C)$ are 
a partition of $\pos(C)$ so the sets $\ri(F)\cap C$ are a partition of $C$.
If $F$ is a face of $\pos(C)$ where $\ri(F)\cap C\neq\emptyset$ then
$\ri(F\cap C)=\ri(F)\cap C$ by (\ref{eq:riint}). This proves that the 
above mapping is a bijection.
\hspace*{\fill}$\Box$\\
\par
The decomposition (\ref{eq:stratum}) justifies a definition:
\begin{Def}
Let $C\subset\bE$ be a convex subset. For every $x\in C$ a unique face 
$F(C,x)$ of $C$ is defined by the condition $x\in\ri(F(C,x))$.
\end{Def}
\par
We describe suprema of faces.
\begin{Lem}
\label{lem:supf}
If $C\subset\bE$ is a convex subset and $\{F_{\alpha}\}_{\alpha\in I}$ 
is a non-empty family of faces of $C$ with $x_{\alpha}\in\ri(F_{\alpha})$ 
for all $\alpha\in I$, then for any 
$z\in\ri(\conv\{x_{\alpha}\mid\alpha\in I\})$ we have
$\bigvee_{\alpha\in I}F_{\alpha}=F(C,z)$.
\end{Lem}
{\em Proof:\/}
Since $z\in F(C,z)$ and since $z$ is in the relative interior of the 
convex set $\conv\{x_{\alpha}\mid\alpha\in I\}$, this convex set is 
included in $F(C,z)$ by (\ref{eq:ri_inc}). So all the $x_{\alpha}$ belong 
to $F(C,z)$. Again by (\ref{eq:ri_inc}) all the faces $F_{\alpha}$ are 
included in $F(C,z)$ because $x_{\alpha}\in\ri(F_{\alpha})$. Thus 
$F(C,z)$ is an upper bound for the family $\{F_{\alpha}\}_{\alpha\in I}$ 
and thus $\bigvee_{\alpha\in I}F_{\alpha}\subset F(C,z)$. Conversely we 
have $z\in\conv\{x_{\alpha}\mid\alpha\in I\}\subset
\bigvee_{\alpha\in I}F_{\alpha}$, so
$F(C,z)\subset\bigvee_{\alpha\in I}F_{\alpha}$ by (\ref{eq:ri_inc}) 
because $z\in\ri(F(C,z))$.
\hspace*{\fill}$\Box$\\
\par
Some faces of $C$ are obtained by intersection of $C$ with a hyperplane, 
these are the {\it exposed faces}. Different to Rockafellar or Schneider 
\cite{rock,schneider} we always include $\emptyset$ and $C$ to the exposed 
faces in order to turn this set into a lattice.
\begin{Def}[Exposed face lattice]
Let $C\subset\bE$ be a convex subset. The {\it support function} of $C$ 
is $\bE\to\bR\cup\{\pm\infty\}$,
$u\mapsto h(C,u):=\sup_{x\in C}\scal{u,x}$.
For non-zero $u\in\bE$
\[
H(C,u):=\{x\in\bE:\scal{u,x}=h(C,u)\}
\]
is an affine hyperplane in $\bE$ unless $H(C,u)=\emptyset$ when
$h(C,u)=-\infty$ with $C=\emptyset$ or $h(C,u)=\infty$, when $C$ is 
unbounded in the direction of $u$. If $H(C,u)\neq\emptyset$, then we call 
it a {\it supporting hyperplane} of $C$. The {\it exposed face} of $C$ by $u$ 
is 
\[
F_{\perp}(C,u):=C\cap H(C,u).
\]
The faces $\emptyset$ and $C$ are exposed faces of $C$ by definition 
called {\it improper exposed faces}. All other exposed faces are {\it 
proper}. The set of exposed faces of $C$ will be denoted by 
$\cF_{\perp}(C)$ called the {\it exposed face lattice} of $C$. A face of 
$C$, which is not an exposed face is a {\it non-exposed face}.
\end{Def}
\par
It is easy to show $\cF_\perp(C)\subset\cF(C)$. An example of a 
non-exposed faces is given in Figure~\ref{fig:exposed_faces_normal_cones}, 
left. It is well-known that the intersection of exposed faces is an 
exposed face, see e.g.\ Schneider \cite{schneider}, but the following 
details were not found in the literature.
\begin{Pro}
\label{pro:ex_lattice}
Let $C\subset\bE$ be a convex set and let $U\subset\bE\setminus\{0\}$ be a 
non-empty set of directions. Then $\ri(\conv(U))\setminus\{0\}$ is 
non-empty and every vector $v$ in this set satisfies 
$\bigcap_{u\in U}F_{\perp}(C,u)=F_{\perp}(C,v)$ unless the intersection is 
empty.
\end{Pro}
{\em Proof:\/} 
Since $U\neq\emptyset$ we have $\ri(U)\neq\emptyset$ (see \cite{rock}, 
Thm.\ 6.2). If we had $\ri(\conv(U))=\{0\}$ then $\conv(U)$ would be 
$\{0\}$, which was excluded in the assumptions. This proves the first 
assertion.
\par
Let $F:=\bigcap_{u\in U}F_{\perp}(C,u)$ and
$G:=\bigcap_{u\in \conv(U)\setminus\{0\}}F_{\perp}(C,u)$.
First we show $F=G$. The non-trivial part is to prove $F\subset G$.
A vector $v\in\conv(U)\setminus\{0\}$ is a convex combination
$v=\sum_{i}\lambda_i u_i$ for $u_i\in U$ and non-negative real
scalars $\lambda_i$ summing up to one. If $x\in F$ then
$x\in F_{\perp}(C,u_i)$ for all $i$ and then
\[\textstyle
\scal{v,x}=\sum_{i}\lambda_i\scal{u_i,x}
=\sum_{i}\lambda_i\max_{s\in C}\scal{u_i,s}
\geq\max_{s\in C}\sum_{i}\lambda_i\scal{u_i,s}
=\max_{s\in C}\scal{v,s},
\] 
so $x\in F_{\perp}(C,v)$. The vector $v$ was arbitrary. So
$x\in G$ and we have $F=G$ indeed.
\par
We assume that $G\neq\emptyset$ and prove $G=F_{\perp}(C,v)$ for 
$v\in\ri(\conv(U))\setminus\{0\}$. To prove the 
non-trivial inclusion $F_{\perp}(C,v)\subset G$ assume by 
contradiction that there is a point $y\in F_{\perp}(C,v)\setminus G$, 
i.e.\ there exists $u_0\in\conv(U)\setminus\{0\}$ such that
\[
y\in F_{\perp}(C,v)\setminus F_{\perp}(C,u_0).
\]
Since $v$ lies in the relative interior of $\conv(U)$ and $u_0$ lies 
in $\conv(U)$ there exists $\lambda\in(0,1)$ and $u_1\in\conv(U)$ such 
that $v=\lambda u_0 + (1-\lambda)u_1$ (see Theorem 6.4 in 
\cite{rock}). We assume that $u_1\neq 0$ by performing a small 
perturbation of this point along the direction $v-u_0$ if necessary. 
Now let $x\in G$. Then we have
$x\in F_{\perp}(C,u_0)\cap F_{\perp}(C,u_1)$ so the estimation
\begin{eqnarray*}
\lefteqn{\textstyle\scal{v,y}
=\lambda\scal{u_0,y}+(1-\lambda)\scal{u_1,y}
< \lambda\max_{z\in C}\scal{u_0,z}+(1-\lambda)\scal{u_1,y}}\\
& \leq & \lambda\scal{u_0,x}+(1-\lambda)\scal{u_1,x} = \scal{v,x}
\hspace{5cm}
\end{eqnarray*}
gives the contradiction $y\not\in F_{\perp}(C,v)$.
\hspace*{\fill}$\Box$\\
\begin{figure}[t!]
\centerline{%
\begin{picture}(4.3,3)
\ifthenelse{\equal{\ftype}{eps}}{%
\put(0,0){\includegraphics*[height=3cm, bb=30 45 280 205]{cake.eps}}}{%
\put(0,0){\includegraphics*[height=3cm, bb=40 55 340 255]{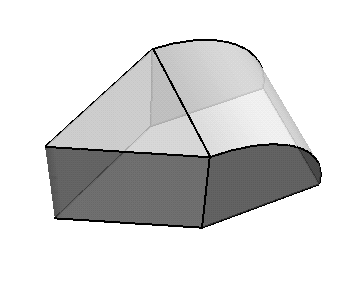}}}
\put(1.6,2.8){$x$}
\put(2.6,0.8){\textcolor{white}{$y$}}
\end{picture}
}
\caption{\label{fig:normal_cones}%
This depicted convex set $K$ is a composition of two right prisms, one 
based on a triangle the other based on a quarter disk. The supremum of the 
extreme points $x$ and $y$ is the the top triangle in $\cF_\perp(K)$ 
and the segment $[x,y]$ in $\cF(K)$.}
\end{figure}
\par
Given a convex subset $C\subset\bE$ the property of a subset of $C$ to be 
an exposed face of $C$ is a closure property by 
Prop.~\ref{pro:ex_lattice}. Thus, by Remark~\ref{rem:closure_complete} the 
exposed face lattice
\begin{equation}
\label{ex_face_lattice}
(\cF_{\perp}(C),\subset,\cap,\vee,\emptyset,C)
\end{equation}
is a complete lattice ordered by inclusion and the infimum is the 
intersection. Although we have the inclusion of 
$\cF_{\perp}(C)\subset\cF(C)$ into the face lattice (\ref{face_lattice}), 
$\cF_{\perp}(C)$ is not in general a sublattice of $\cF(C)$. Both lattices 
have the intersection as infimum but their suprema may be different. An 
example is drawn in Figure~\ref{fig:normal_cones}.
\par
We prove a technical detail for the next assertion. If $C$ is convex 
subset of $\bE$, $x\in\bE$ and $\{x\}\subsetneq C$ then the equality
\begin{equation}
\label{eq:too_many_details}
\ri(\conv(C\setminus\{x\}))=\ri(C)
\end{equation}
holds. If $C\setminus\{x\}$ is not convex then 
$\conv(C\setminus\{x\})=C$ and the equality follows. If 
$C\setminus\{x\}$ is convex then $x$ is an extreme point of $C$. Hence,
unless $C=\{x\}$, we have $\ri(C)\subset C\setminus\{x\}\subset C$. 
Therefore $C\setminus\{x\}$ lies between $\ri(C)$ and the closure 
$\overline{C}$ of $C$. Thus the relative interiors of $C\setminus\{x\}$ 
and $C$ are equal by Coro.\ 6.3.1 in \cite{rock}.
\begin{Cor}
\label{cor:ex_inf}
Let $C,D\subset\bE$ be convex subsets. If $D$ contains a non-zero vector, 
then $\ri(D)$ contains a non-zero vector. If 
$\bigcap_{u\in D\setminus\{0\}}F_{\perp}(C,u)\neq\emptyset$ then this 
intersection is the exposed face $F_{\perp}(C,v)$ for any non-zero 
$v\in\ri(D)$.
\end{Cor}
{\em Proof:\/}
By Prop.~\ref{pro:ex_lattice} we have for any vector
$v\in\ri(\conv(D\setminus\{0\}))\setminus\{0\}$ the equality of the 
intersection with the face $F_{\perp}(C,v)$. With 
(\ref{eq:too_many_details}) applied to $x:=0$ and $C:=D$ we get 
$\ri(\conv(D\setminus\{0\}))\setminus\{0\}=\ri(D)\setminus\{0\}$.
\hspace*{\fill}$\Box$\\
%
%
%
%
%
\section{Normal cones}
\label{sec:normal_cones}
\par
We study normal cones of a convex subset $C\subset\bE$ of the 
finite-dimensional real Euclidean vector space 
$(\bE,\langle\cdot,\cdot\rangle)$. There is an antitone lattice 
isomorphism between exposed faces and normal cones. 
\begin{Def}
The {\it normal cone} of $C$ at $x\in C$ is 
\begin{equation}
\label{normal_cone_definition}
\nc(C,x):=\{u\in\bE:\scal{u,y-x}\leq 0\text{ for all }y\in C\}
\end{equation}
and vectors in $\nc(C,x)$ are called {\it normal vectors} of $C$ at $x$. 
\end{Def}
\par
There is a pointwise relation between exposed faces and normal cones. If 
$C\subset\bE$ is a convex subset, then for arbitrary $x\in C$ and non-zero 
$u\in\bE$ the equivalence of the following statements is easy to prove.
\begin{equation}
\label{eq:ptw_duality}
\begin{array}{crcl}
\bullet & \scal{u,x} & = & h(C,u),\\
\bullet & x & \in & F_{\perp}(C,u),\\
\bullet & u & \in & \nc(C,x).
\end{array}
\end{equation}
The following relations are easy to prove by elementary means. If 
$F\subset C$ is a convex subset, $x\in\ri(F)$ and $y\in C$, then we have
\begin{equation}
\label{eq:nc_for_begn}
\begin{array}{rl}
\text{(i)} & \nc(C,x)\perp\lin(F),\\
\text{(ii)} & \text{if } y\in F \text{ then }  
\nc(C,y)\supset\nc(C,x),\\
\text{(iii)} & \text{if } y\in\ri(F) \text{ then } 
\nc(C,y)=\nc(C,x),\\
\text{(iv)} & \text{if } u,-u\in\nc(C,y) \text{ then } 
u\in\lin(C)^{\perp}.
\end{array}
\end{equation}
The orthogonal complement with respect to the Euclidean inner product is 
denoted by ${}^{\perp}$.
\begin{Lem}
\label{lin_complement}
Let $x\in C$. Then $\nc(C,x)=(\nc(C,x)\cap\lin(C))+\lin(C)^{\perp}$ holds 
and the following statements are equivalent.\\[0.5cm]
\centerline{\begin{tabular}{cl}
$\bullet$ & the normal cone $\nc(C,x)$ is a vector space,\\
$\bullet$ & $x\in\ri(C)$,\\
$\bullet$ & $\nc(C,x)=\lin(C)^{\perp}$.
\end{tabular}}
\end{Lem}
{\it Proof:\/} 
Let $x\in C$. The direct sum decomposition of $\nc(C,x)$ follows from
$\nc(C,x)+\lin(C)^\perp\subset\nc(C,x)$. Since $\nc(C,x)$ is a convex 
cone, it is sufficient to prove the inclusion 
$\lin(C)^{\perp}\subset\nc(C,x)$: if $u\in\lin(C)^{\perp}$ 
then $\scal{u,y-x}=0$ for all $y\in C$ so $u\in\nc(C,x)$.
\par
Now let us assume that $\nc(C,x)$ is a vector space. Then for 
$u\in\nc(C,x)$ we have $\pm u\in\nc(C,x)$ and by (\ref{eq:ptw_duality}) we 
get
\[
h(C,u)=\scal{u,x}=-\scal{-u,x}=-h(C,-u).
\]
Thus, for the vectors $u\in\bE$ with $h(C,u)\neq -h(C,-u)$ follows 
$u\not\in\nc(C,x)$, which means $\scal{u,x}<h(C,u)$ by 
(\ref{eq:ptw_duality}). These are exactly the assumption of Theorem 13.1 
in \cite{rock} to prove that $x\in\ri(C)$. Clearly, if $x\in\ri(C)$ then 
$\nc(C,x)=\lin(C)^{\perp}$.
\hspace*{\fill}$\Box$\\
\begin{Def}
\label{def:normal_cone_of_face}
The {\it normal cone} of a non-empty convex subset $F$ of $C$ is defined as
\begin{equation}
\label{def:nml_cone_f}
\nc(C,F):=\nc(C,x)
\end{equation}
for any $x\in\ri(F)$. This definition is consistent by (iii) in
(\ref{eq:nc_for_begn}). The {\it normal cone} of the empty set is defined 
as the ambient space $\nc(C,\emptyset):=\bE$. The 
{\it normal cone lattice} of $C$ is the set of normal cones of all faces 
$\cN(C):=\{\nc(C,F)\mid F\in\cF(C)\}$. We consider the normal cone 
lattice as a poset ordered by set inclusion. The cones $\lin(C)^\perp$ and 
$\bE$ are the {\it improper} normal cones, all other normal cones are
{\it proper}.
\end{Def}
\begin{figure}[t!]
\centerline{%
\begin{picture}(4.2,2.1)
\ifthenelse{\equal{\ftype}{eps}}{%
\put(0.3,0){\includegraphics[height=2.0cm,bb=0 0 300 156]{%
B040non_exposed_faces.eps}}}{%
\put(0.3,0){\includegraphics[height=2.0cm,bb=0 0 300 185]{%
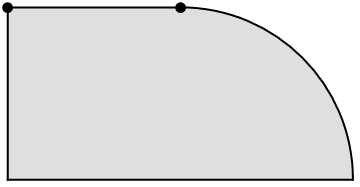}}}
\put(0,1.8){$x$}
\put(2.1,1.6){$y$}
\put(1,0.5){$C$}
\end{picture}}
\caption[Normal cones for non-exposed faces]{%
\label{fig:normal_cone_counter}
The union $C$ of a square and a quarter disk with extreme points $x,y$:
$\{x\}$ is an exposed face while $\{y\}$ is a non-exposed face.  The face 
$\{y\}$ has the same normal cone as the face $[x,y]$. The normal cone of 
$\{y\}$ is included in the normal cone of $\{x\}$, even though $\{x\}$ and 
$\{y\}$ are unrelated in the partial ordering of inclusion.}
\end{figure}
\par
The assignment of normal cones to faces $\cF(C)\to\cN(C)$, 
$F\mapsto\nc(C,F)$ is an antitone mapping between posets. This follows 
from (ii) in (\ref{eq:nc_for_begn}). But the faces of two included normal 
cones may be unrelated, see Figure~\ref{fig:normal_cone_counter}. We work 
towards the antitone lattice isomorphism $\cF_\perp(C)\to\cN(C)$. 
\begin{Lem}
\label{lem:face_cone_dual_u}
If $F\in\cF(C)$ is a face and $u\in\bE\setminus\{0\}$ then
$F\subset F_{\perp}(C,u)$ if and only if  $u\in\nc(C,F)$. For all 
$u\in\bE\setminus\{0\}$ we have $u\in\nc(C,F_{\perp}(C,u))$.
\end{Lem}
{\em Proof:\/}
The assertion is trivial for $F=\emptyset$. Otherwise let us assume that 
the inclusion $F\subset F_{\perp}(C,u)$ holds and consider a point 
$x\in\ri(F_{\perp}(C,u))$. We have $u\in\nc(C,x)=\nc(F_{\perp}(C,u))$ by 
the relation (\ref{eq:ptw_duality}) and by definition 
(\ref{def:nml_cone_f}) of a normal cone. Since $F\subset F_{\perp}(C,u)$
we have $\nc(C,F_{\perp}(C,u))\subset\nc(C,F)$ by the antitone normal 
cone assignment. Conversely, if $u\in\nc(C,F)$ then 
for $x\in\ri(F)$ we have $u\in\nc(C,x)$. Thus $x\in F_{\perp}(C,u)$ by 
the relation (\ref{eq:ptw_duality}) and 
(\ref{eq:ri_inc}) gives $F\subset F_{\perp}(C,u)$. The second assertion is 
the special case of $F=F_\perp(C,u)$.
\hspace*{\fill}$\Box$\\
\par
We consider smallest upper bounds of exposed faces for arbitrary subsets 
of $C$. This is consistent by completeness (\ref{ex_face_lattice}) of the 
exposed face lattice $\cF_{\perp}(C)$:
\begin{Def}
The {\it smallest exposed face} of $C$ that contains a subset $F\subset C$ 
is
\begin{equation}
\label{eq:smallest_exp_face}
\sex(F):=\bigcap\{G\in\cF_{\perp}(C)\mid F\subset G\}.
\end{equation}
\end{Def}
\par
Properties of the smallest exposed face are:
\begin{Lem}
\label{lem:sup_exp}
If $F\in\cF(C)$ is a proper face, then 
$\sex(F)=\bigcap_{u\in\nc(C,F)\setminus\{0\}}F_{\perp}(C,u)$ is a proper 
exposed face. We have $\ri(\nc(C,F))\neq\{0\}$ and for
each non-zero $v\in\ri(\nc(C,F))$ we have $\sex(F)=F_{\perp}(C,v)$.
If $F\in\cF(C)$ is a face then $\nc(C,\sex(F))=\nc(C,F)$. 
\end{Lem}
{\em Proof:\/}
By Lemma~\ref{lem:face_cone_dual_u} if $u\in\bE$ is non-zero, then the 
face $F$ is included in $F(C,u)$ if and only if $u\in\nc(C,F)$.
\par
Relative interior points of the proper face $F$ do not belong to $\ri(C)$,
so by Lemma~\ref{lin_complement} the normal cone of $F$ is strictly larger 
than $\lin(C)^\perp=\nc(C,C)$. Choosing any
$u\in\nc(C,F)\setminus\lin(C)^\perp$ we get that $F$ but not $C$ is 
included in $F_\perp(C,u)$. So $\sex(F)$ is a proper exposed face of $C$
and the intersection expression for $\sex(F)$ follows. As 
$F\neq\emptyset$, any non-zero vector $v\in\ri(\nc(C,F))$ gives 
$\sex(F)=F_{\perp}(C,v)$ by Cor.~\ref{cor:ex_inf}.
\par
Since $F\subset\sex(F)$, the inclusion $\nc(C,\sex(F))\subset\nc(C,F)$ 
follows from antitone assignment of normal cones. For every non-zero 
vector $u\in\nc(C,F)$ we have $F\subset F_{\perp}(C,u)$. Hence
$\sex(F)\subset F_{\perp}(C,u)$ and so $u\in\nc(C,\sex(F))$.
\hspace*{\fill}$\Box$\\
\par
We arrive at the main results of this section.
\begin{Pro}
\label{lattice_antitone_iso}
Assume that $C$ has not exactly one point. Then the assignment of 
normal cones to exposed faces $\nc(C):\cF_{\perp}(C)\to\cN(C)$, 
$F\mapsto \nc(C,F)$ is  an antitone lattice isomorphism.
\end{Pro}
{\em Proof:\/}
The two lattices $\cF_{\perp}(C)$ and $\cN(C)$ are partially ordered by 
set inclusion. They are linked by the antitone mapping of posets
\[
\nc(C)|_{\cF_{\perp}(C)}:\cF_{\perp}(C)\to\cN(C),
\quad F\mapsto \nc(C,F).
\]
This mapping is surjective because
a face $F$ of $C$ has the same normal cone as the smallest exposed 
face that contains $F$, see Lemma~\ref{lem:sup_exp}.
\par
We can show that $\nc(C)|_{\cF_{\perp}(C)}$ has an antitone inverse. Then 
Remark~\ref{rem:isotone_isomorphism} implies that 
$\nc(C)|_{\cF_{\perp}(C)}$ is an (antitone) lattice isomorphism. Let us 
prove that this map is injective and consider two proper exposed
faces $F,G$ of $C$ with the same normal cone $N$. Then there 
exists by Lemma~\ref{lin_complement} a non-zero vector $u\in N$, so there 
is a non-zero $v\in\ri(N)$. As $F,G\neq\emptyset$, 
Lemma~\ref{lem:sup_exp} proves that $F=F_{\perp}(C,v)=G$. By  
Lemma~\ref{lin_complement} only the improper face $C$ has the smallest 
possible normal cone $\lin(C)^{\perp}$. It remains to show that 
$\nc(C,F)=\bE$ implies $F=\emptyset$ for an exposed face $F$ of $C$. If 
$\nc(C,F)=\bE$ holds for a non-empty face $F$ then 
Lemma~\ref{lin_complement} shows that $F=C$ and
$\lin(C)=\bE^{\perp}=\{0\}$. Thus, $C$ has exactly one point but this 
case was excluded in the assumptions.
\par
We show that the inverse $\cN(C)\to\cF_\perp(C)$ is antitone. For proper 
exposed faces $F,G$ of $C$ the inclusion  $\nc(G)\subset\nc(F)$ implies 
$\sex(F)\subset\sex(G)$ by Lemma~\ref{lem:sup_exp}. As $F,G$ are exposed 
we have $F=\sex(F)$ and $G=\sex(G)$, hence $F\subset G$. The greatest 
element $\bE$ of $\cN(C)$ maps to the smallest element $\emptyset$ of 
$\cF_\perp(C)$ and the smallest element $\lin(C)^{\perp}$ of $\cN(C)$ maps 
to the greatest element $C$ of $\cF_\perp(C)$.
\hspace*{\fill}$\Box$\\
\par
By definition of the normal cone of a face and by antitone assignment of 
normal cones the isomorphism $\cF_\perp(C)\to\cN(C)$ in 
Prop.~\ref{lattice_antitone_iso} is for proper exposed faces 
$F\in\cF_\perp(C)$
\begin{equation}
\label{eq:anti_isos}
\begin{array}{rcl}
F & \mapsto & \bigcap_{x\in F}\nc(C,x)
\quad = \quad \nc(C,y) \qquad\text{for any }y\in\ri(F),\\
N & \mapsto & \bigcap_{u\in N\setminus\{0\}}F_\perp(C,u)
\quad = \quad F_\perp(C,v) \qquad\text{for any }v\in\ri(N)\setminus\{0\}.
\end{array}
\end{equation}
The second mapping defined for proper normal cones $N\in\cN(C)$ describes 
the inverse $\cN(C)\to\cF_\perp(C)$ by Lemma~\ref{lem:sup_exp}.
Now we shows that intersections of normal cones are normal cones, so by 
Remark~\ref{rem:closure_complete} the normal cone lattice is a complete 
lattice with intersection as the infimum
\begin{equation}
\label{eq:normal_l}
(\cN(C),\subset,\cap,\vee,\lin(C)^\perp,\bE).
\end{equation}
\begin{Pro}
\label{normal_cone_intersection}
If $\{N_\alpha\}_{\alpha\in I}\subset\cN(C)$ is a non-empty family of 
normal cones, then
$\bigwedge_{\alpha\in I}N_\alpha=\bigcap_{\alpha\in I}N_\alpha$ and this 
intersection is a face of $N_{\widetilde{\alpha}}$ for every 
$\widetilde{\alpha}\in I$ with $N_{\widetilde{\alpha}}\neq\bE$.
\end{Pro}
{\em Proof:\/}
As $\bE$ is the greatest element of $\cN(C)$ we assume $N_\alpha\neq\bE$ 
for all $\alpha\in I$ and we assume that $C$ has not exactly one point,
without restricting generality. As $\nc(C,\emptyset)=\bE$ we 
choose throughout for $\alpha\in I$ a family of (non-empty) faces 
$F_\alpha$ with $\nc(C,F_\alpha)=N_\alpha$. Let 
$x_\alpha\in\ri(F_\alpha)$ for $\alpha\in I$ and let 
$z\in\ri(\conv\{x_{\alpha}\mid\alpha\in I\})$. So Lemma~\ref{lem:supf} 
shows $F(C,z)=\bigvee_{\alpha\in I}F_{\alpha}$. By 
Prop.~\ref{lattice_antitone_iso} we have 
$K:=\bigwedge_{\alpha\in I}\nc(C,F_{\alpha})
=\nc(C,\bigvee_{\alpha\in I}F_{\alpha})=\nc(C,z)$.
\par
The assignment of a normal cone is antitone, so for all 
$\widetilde{\alpha}\in I$ we have $K\subset\nc(C,F_{\widetilde{\alpha}})$. 
This proves one inclusion, it remains to show
$\bigcap_{\alpha\in I}\nc(C,F_{\alpha})\subset K$. We write $z$ as a 
convex combination for $n\in\bN$, $\lambda_i>0$ and $\alpha(i)\in I$ for 
$i=1,\ldots,n$ in the form $z=\sum_{i=1}^n\lambda_ix_{\alpha(i)}$.
Hence, if $u\in\bigcap_{\alpha\in I}\nc(C,F_{\alpha})$, then for all
$x\in C$ we have the inequality
$\langle u,x-z\rangle=\sum_{i=1}^n\lambda_i\langle u, x-x_{\alpha(i)}
\rangle\leq 0$. This proves $u\in\nc(C,z)$.
\par
For $\widetilde{\alpha}\in I$ let us prove that $K$ is a face of 
$N_{\widetilde{\alpha}}$. Let $u,v,w\in N_{\widetilde{\alpha}}$, $v\in K$ 
and $v\in\,]u,w[$. If $u=0$ then $w=\lambda v$ for some $\lambda>0$, then 
$u,w\in K$ because $K$ is a convex cone including $v$. If  $u,w\neq 0$ and 
$v=0$ then $u,w\in\lin(C)^{\perp}$. By Lemma~\ref{lin_complement} the 
vector space $\lin(C)^{\perp}$ belongs to every normal cone of $C$, so 
$u,w\in K$. Let us assume $u,v,w\neq 0$. For every $\alpha\in I$ holds
$v\in N_\alpha=\nc(C,F_\alpha)$ so $F_\alpha\subset  F_{\perp}(C,v)$ by 
Lemma~\ref{lem:face_cone_dual_u}. Now Prop.~\ref{pro:ex_lattice} shows 
$F_{\perp}(C,v)=F_{\perp}(C,u)\cap F_{\perp}(C,w)$, so we have
\[
F_\alpha\subset F_{\perp}(C,v)=F_{\perp}(C,u)\cap F_{\perp}(C,w)
\subset F_{\perp}(C,u).
\]
This gives $\nc(C,F_{\perp}(C,u))\subset\nc(C,F_\alpha)$ and 
Lemma~\ref{lem:face_cone_dual_u} completes the proof
with $u\in\nc(C,F_{\perp}(C,u))$. The proof of $w\in\nc(C,F_\alpha)$ is a 
complete analogue.
\hspace*{\fill}$\Box$\\
%
%
%
%
%
%
%
\section{Cylinders}
\label{sec:cylinders}
\par
This section explains a lifting construction for projections of 
convex sets. Lifting is an isomorphism for face lattices, we characterize 
lifted faces. As an application, the projections of state spaces 
introduced in Section~\ref{sec:context} are decomposed by Weis \cite{Weis} 
using this lifting. Throughout this section let $C$ be a convex 
subset of a finite-dimensional real Euclidean vector space 
$(\bE,\scal{\cdot,\cdot})$ and let $V$ be a linear subspace of $\bE$.
\par
If $\emptyset\neq\bA\subset\bE$ is an affine subspace, then there is a 
unique affine mapping $\pi_{\bA}:\bE\rightarrow\bA$, called the
{\it orthogonal projection} to $\bA$, such that for all $x\in\bE$  we have
\begin{equation}
\label{orthogonal_projection_definition}
(x-\pi_{\bA}(x))\perp\lin(\bA).
\end{equation}
\par
We study the orthogonal projection $\pi_V:\bE\to V$ to $V$. This, thought 
of as acting on sets, may be written for $M\subset\bE$ in the form
\begin{equation}
\label{eq:projection_cylinder}
\pi_{V}(M)=(M+V^{\perp})\cap V.
\end{equation}
In addition to the projection 
$\pi_V(C)$ we will study the cylinder $C+V^{\perp}$, which connects 
the projection $\pi_V(C)$ to $C$.
\par
There is a basic tool for the study of cylinders, which is reminiscent
of the modular law for lattices (\ref{eq:mod_l}).
\begin{Lem}
\label{cyl:modular}
Let $X,Y,Z\subset\bE$ such $Z\pm X\subset Z$. Then
$X+(Y\cap Z)=(X+Y)\cap Z$.
\end{Lem}
{\em Proof:\/}
The inclusion $(X+Y)\cap Z\subset X+(Y\cap Z)$ is proved by taking 
vectors $x\in X$ and $y\in Y$ such that $x+y\in Z$. Then
$y=(x+y)-x\in Z$. For the converse $X+(Y\cap Z)\subset(X+Y)\cap Z$ we 
choose vectors $x\in X$ and $t\in Y\cap Z$. Then $t+x\in Z$.
\hspace*{\fill}$\Box$\\
\par
A special case of Lemma~\ref{cyl:modular} is the
{\it modular law for affine spaces}. Let $\bA\subset\bE$ be an affine 
subspace with translation vector space $\lin(\bA)$. If 
$X\subset\lin(\bA)$ then for arbitrary $Y\subset\bE$ we have
\begin{equation}
\label{eq:mod_a}
X+(Y\cap\bA)=(X+Y)\cap\bA.
\end{equation}
We will use this modular law as indicated in Figure~\ref{fig:modular}.
\begin{figure}[t!]
\centerline{%
\begin{picture}(5.8,5)
\ifthenelse{\equal{\ftype}{eps}}{%
\put(0,0){\includegraphics[height=5cm, bb=15 130 390 460, 
clip=]{modular.eps}}}{%
\put(0,0){\includegraphics[height=5cm, bb=20 160 500 570, 
clip=]{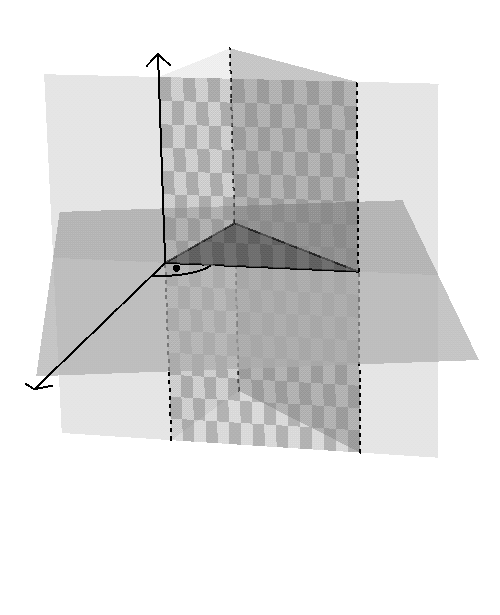}}}
\put(2.5,2.4){\textcolor{white}{$C$}}
\put(4.8,1.2){$V$}
\put(-0.2,0.8){$v$}
\put(4.6,3.9){$H$}
\put(1.0,4.7){$V^\perp$}
\end{picture}}
\caption{\label{fig:modular}We start with a plane $V$ and an arbitrary 
subset $C$ in $\bR^3$. For simplicity in the drawing we choose $C$ a 
triangle in $V$. A non-zero vector $v\in V$ defines the supporting 
hyperplane $H=H(C,v)$ with $v\perp H$. We have 
$V^{\perp}\subset\{v\}^{\perp}=\lin(H)$. So by the modular law for affine 
spaces $V^{\perp}+(C\cap H)=(V^{\perp}+C)\cap H$ holds. This set is drawn 
tiled.}
\end{figure}
\begin{Def}
We define the {\it lift} from $V$ to $C$ (or along $V^{\perp}$ to $C$)
as the mapping $L_V^C: 2^{\bE}\to 2^C$, $M \mapsto (M+V^{\perp})\cap C$.
Here $2^{\bE}$ denotes the power set of $\bE$ and $2^C$ the power set of 
$C$.
\end{Def}
\begin{Lem}
\label{ap_lift_iso}
The projection $\pi_V:2^{\bE}\to 2^V$ is isotone with respect to set 
inclusion and we have
\[
L_V^C=L_V^C\circ L_V^C=L_V^C\circ\pi_V.
\]
If $\cM$ is a family of subsets
of $\pi_V(C)$, then $\pi_V$ is left inverse to $L_V^C|_{\cM}$.
In particular
\[
L_V^C|_{\cM}:\cM\to\{L_V^C(M):M\in\cM\}
\]
is a bijection. The mapping $L_V^C|_{\cM}$ is an isotone isomorphism of 
posets (partially ordered by set inclusion). 
\end{Lem}
{\em Proof:\/}
Trivial.\hspace*{\fill}$\Box$\\
\begin{Lem}[Lifted faces]
\label{lemma:lifted_faces}
If $F$ is a face of $\pi_V(C)$ then the lift $L_V^C(F)$ is a face of 
$C$. The exposed face for non-zero $v\in V$ transforms according to 
$L_V^C(F_{\perp}(\pi_V(C),v))=F_{\perp}(C,v)$.
\end{Lem}
{\em Proof:\/}
For a face $F$ of $\pi_V(C)$ we show that $L_V^C(F)$ is a face of $C$. 
To this aim we choose $x,y,z\in C$ such that $y\in]x,z[$ and 
$y\in L_V^C(F)$. We have to prove $x,z\in L_V^C(F)$. By (\ref{eq:cori})
the projection $\pi_V$ commutes with reduction to the relative interior 
of a convex set, so we have $\pi_V(y)\in]\pi_V(x),\pi_V(z)[$. Since
$y\in L_V^C(F)$ we have $\pi_V(y)\in F$. Since $F$ is a face we obtain
$\pi_V(x),\pi_V(z)\in F$. Then
\[
x\in L_V^C\circ\pi_V(x)
=(\pi_V(x)+V^{\perp})\cap C\subset(F+V^{\perp})\cap C
=L_V^C(F).
\]
Analogously we have $z\in L_V^C(F)$, so $L_V^C(F)$ is a face of $C$.
\par
The support functions of $C$ and $\pi_V(C)$ are equal on $V$ because for 
all $x\in\bE$ and $v\in V$ we have $\scal{v,x}=\scal{v,\pi_V(x)}$. If 
$v\in V$ is a non-zero vector then the hyperplanes $H(C,v)$ and 
$H(\pi_V(C),v)$ are equal. Since $v\in V$ we have 
$V^{\perp}\subset\{v\}^{\perp}=\lin(H(\pi_V(C),v))$ and we can 
apply the modular law for affine spaces (\ref{eq:mod_a}) as follows 
\begin{eqnarray*}
\lefteqn{%
V^{\perp}+F_{\perp}(\pi_V(C),v)
=V^{\perp}+[\pi_V(C)\cap H(\pi_V(C),v)]
=[V^{\perp}+\pi_V(C)]\cap H(\pi_V(C),v)}\\
& = & (V^{\perp}+C)\cap H(C,v).\hspace{10cm}
\end{eqnarray*}
This gives
\begin{eqnarray*}
\lefteqn{
L_V^C(F_{\perp}(\pi_V(C),v))=(F_{\perp}(\pi_V(C),v)+V^{\perp})\cap 
= (V^{\perp}+C) \cap H(C,v) \cap C}\\
& = & H(C,v)\cap C = F_{\perp}(C,v)\hspace{8cm}
\end{eqnarray*}
finally.
\hspace*{\fill}$\Box$\\
\begin{Def}
With respect to $C$ and $V$, the face $L_V^C(F)\in\cF(C)$
is called the {\it lifted face} of $F\in\cF(\pi_V(C))$.
The {\it lifted face lattice} is
\[
\cF_V^C:=\{L_V^C(F):F\in\cF(\pi_V(C))\}.
\]
The {\it lifted exposed face lattice} is
\begin{equation}
\label{cyl:lifted_exposed_lattice}
\cF_{V,\perp}^C:=\{L_V^C(F):F\in\cF_{\perp}(\pi_V(C))\}
\end{equation}
where $\cF(\pi_V(C))$ is the face lattice of $\pi_V(C)$ and
$\cF_{\perp}(\pi_V(C))$ is the exposed face lattice of
$\pi_V(C)$. We consider $\cF_V^C$ and $\cF_{V,\perp}^C$
partially ordered by set inclusion. 
\end{Def}
\begin{Pro}[Lifted face lattices]
\label{pro:ll_iso}
The lifts from $V$ to $C$ restricted to the face lattices of $\pi_V(C)$ are
\[\begin{array}{rrcl}
L_V^C|_{\cF(\pi_V(C))}: & \cF(\pi_V(C)) &
\to & \cF_V^C\subset\cF(C),\\
L_V^C|_{\cF_{\perp}(\pi_V(C))}: & \cF_{\perp}(\pi_V(C)) &
\to & \cF_{V,\perp}^C\subset\cF_{\perp}(C).
\end{array}\]
These mappings are isotone lattice isomorphisms. The infimum in the lifted 
face lattices is the intersection.
\end{Pro}
{\em Proof:\/}
The mapping $L_V^C$ restricted to $\cF(\pi_V(C))$ resp.\ to 
$\cF_{\perp}(\pi_V(C))$ is a bijection to $\cF_V^C$ resp.\  to 
$\cF_{V,\perp}^C$ by Lemma~\ref{ap_lift_iso}. The ranges are included in 
the face lattice of $C$ resp.\ in the exposed face lattice of $C$ 
by Lemma~\ref{lemma:lifted_faces}.
\par
The mappings $L_V^C$ and $\pi_V$ (on the considered domains) are inverse 
to each other and they are isotone  with respect to set inclusion by 
Lemma~\ref{ap_lift_iso}. Hence the lift is a lattice isomorphism in 
each case by Remark~\ref{rem:isotone_isomorphism}. 
\par
Finally, by direct sum decomposition of $\bE=V+V^{\perp}$ we have for a 
non-empty family $\{F_{\alpha}\}_{\alpha\in I}$ of faces of $\pi_V(C)$
\[\textstyle
L_V^C(\bigcap_{\alpha\in I}F_{\alpha})
=(\bigcap_{\alpha\in I}F_{\alpha}+V^{\perp})\cap C
=\bigcap_{\alpha\in I}(F_{\alpha}+V^{\perp})\cap C
=\bigcap_{\alpha\in I}L_V^C(F_{\alpha}),
\]
the infimum in the lifted face lattices is the intersection.
\hspace*{\fill}$\Box$\\
\par
We notice that the lifted exposed face lattice $\cF_{V,\perp}^C$ is not a 
sublattice of the face lattice $\cF(C)$ because the supremum of lifted 
faces in $\cF(C)$ is not necessarily a lifted face. An example is a 
triangle projected to the linear span of one of its sides, say $c$. Then 
the corners $A$ and $B$ of $c$ belong to $\cF_{V,\perp}^C$, but $c$ does 
not. We characterize the lifted face lattice:
\begin{Pro}[Lift invariant faces]
\label{prop:lift_face_char}
A face $F\in\cF(C)$ belongs to the lifted face lattice $\cF_V^C$ if and 
only if $L_V^C(F)=F$. 
\end{Pro}
{\em Proof:\/}
Let us choose a face $F\in\cF(C)$. If $F$ belongs to $\cF_V^C$ 
then there is a face $G\in\cF(\pi_V(C))$ such that $F=L_V^C(G)$. With 
Lemma~\ref{ap_lift_iso} we obtain 
\[
L_V^C(F)=L_V^C\circ L_V^C(G)=L_V^C(G)=F.
\]
\par
For the converse we assume that $F=L_V^C(F)$. If $\pi_V(F)$ is a face of 
$\pi_V(C)$ then we have $F=L_V^C\circ\pi_V(F)$ and so $F$ is a lifted face.
It remains to prove $\pi_V(F)\in\cF(\pi_V(C))$. To this end let 
$x,y,z\in\pi_V(C)$ such that $y\in]x,z[$ and
$y\in\pi_V(F)$. We must show $x,z\in\pi_V(F)$. We choose 
$\widetilde{x}\in L_V^C(x)$ and $\widetilde{z}\in L_V^C(z)$. Then
$[\widetilde{x},\widetilde{z}]\stackrel{\pi_V}{\longrightarrow}[x,z]$ is 
a bijection so there exists
$\widetilde{y}\in]\widetilde{x},\widetilde{z}[\cap L_V^C(y)$. Since 
$y\in\pi_V(F)$ we have 
$\widetilde{y}\in L_V^C\circ\pi_V(F)=L_V^C(F)=F$ and this proves
$\widetilde{x},\widetilde{z}\in F$ because $F$ is a face of $C$. Then
$x=\pi_V(\widetilde{x})$ and $z=\pi_V(\widetilde{z})$ belong to 
$\pi_V(F)$ and we have proved that $\pi_V(F)$ is a face of $\pi_V(C)$.
\hspace*{\fill}$\Box$\\
\par
There is a canonical space to project onto.
\begin{Cor}
Let $U$ be the orthogonal projection of $V$ onto the the vector space of 
$C$, i.e.\ $U:=\pi_{\lin(C)}(V)$. Then for all $F\subset C$ we have
$L_V^C(F)=L_U^C(F)$. In particular $\cF_V^C=\cF_U^C$ holds.
\end{Cor}
{\em Proof:\/}
We put $W:=\lin(C)$. By straight forward calculation we have for any 
$F\subset C$
\[\textstyle
L_U^C(F)=((V^{\perp}\cap W)+(F+W^{\perp}))\cap\aff(C)\cap C.
\]
By the modular law (\ref{eq:mod_a}) applied to the first two intersection
sets this simplifies to $L_V^C(F)$. The second statement follows now from
Prop.~\ref{prop:lift_face_char}.
\hspace*{\fill}$\Box$\\
\par
Finally we write down the normal cones.
\begin{Lem}[Normal cones]
\label{lemma:normal_cone_cylinder}
Let $a\in C+V^{\perp}$. Then 
$\nc(\pi_V(C),\pi_V(a))=\nc(C+V^{\perp},a)+V^{\perp}$. If $a$ belongs 
to $C$ then $\nc(C+V^{\perp},a)=\nc(C,a)\cap V$. 
\end{Lem}
{\em Proof:\/}
Let $a\in C+V^{\perp}$. We use the relation (\ref{eq:ptw_duality}) to 
prove the first identity. We decompose a vector $u\in\bE$ in the form 
$u=v+w\in\bE$ for $v\in V$ and $w\in V^{\perp}$. If 
$u\in\nc(\pi_V(C),\pi_V(a))$ then 
\[
h(C+V^{\perp},v)=h(\pi_V(C),v)=h(\pi_V(C),u)=\scal{u,\pi_V(a)}
=\scal{v,\pi_V(a)}=\scal{v,a},
\]
so $v\in\nc(C+V^{\perp},a)$ and $u\in\nc(C+V^{\perp},a)+V^{\perp}$. 
Conversely, if $v\in\nc(C+V^{\perp},a)$ then
\[
\scal{u,\pi_V(a)}=\scal{v,\pi_V(a)}=\scal{v,a}
=h(C+V^{\perp},v)=h(\pi_V(C),v)=h(\pi_V(C),u),
\]
so $u\in\nc(\pi_V(C),\pi_V(a))$.
\par
The second equation is as follows. If $u\in\nc(C+V^{\perp},a)$, then 
$u\in\nc(C,a)$ because there are less conditions on normal cones for the 
smaller set $C$. For all $w\in V^{\perp}$ we have
$\langle u,\pm w\rangle\leq 0$ so $u\in V$. Conversely, if 
$u\in\nc(C,a)\cap V$, then for all $x\in C$ and for all $w\in V^{\perp}$ 
we have $\langle u,x+w-a\rangle=\langle u,x-a\rangle\leq 0$ and this 
proves $u\in\nc(C+V^{\perp},a)$.
\hspace*{\fill}$\Box$\\
%
%
%
%
%
\section{Sharp relations}
\label{sec:sharp}
\par
Let $(\bE,\langle\cdot,\cdot\rangle)$ be a 
finite-dimensional real Euclidean vector space and $C\subset\bE$ a convex 
subset. There is a relation (\ref{eq:ptw_duality}) between exposed faces 
and normal cones, this is for $x\in C$ and $u\in\bE\setminus\{0\}$ 
\[
x\in F_{\perp}(C,u)\quad\iff\quad u\in\nc(C,x).
\]
We define two alterations: 
\begin{Def}
A vector $u\in\bE\setminus\{0\}$ is {\it sharp normal} for $C$ if
\begin{equation}
\label{eq:sharp_duality_1}
x\in\ri(F_{\perp}(C,u))\quad\implies\quad
u\in\ri(\nc(C,x)).
\end{equation}
A point $x\in C$ is {\it sharp exposed} in $C$ if
\begin{equation}
\label{eq:sharp_duality_2}
u\in\ri(\nc(C,x))\setminus\{0\}\quad\implies\quad
x\in\ri(F_{\perp}(C,u)).
\end{equation}
\end{Def}
\par
A connection of sharp normal vectors to normal cones will be shown in the 
following section. In this section we show that the above definitions do 
not depend on the ambient space (through the normal cones), the argument 
for sharp exposed points connects these to exposed faces. We show that 
sharp normal vectors are preserved under orthogonal projection of a convex 
set and sharp exposed points are preserved under intersection. An example 
where both (\ref{eq:sharp_duality_1}) and (\ref{eq:sharp_duality_2}) hold 
is a state space:
\begin{Exa}
\label{ex:state_spaces}
For $n\in\bN$ let ${\rm Mat}(\bC,n)$ be the set of complex $n\times n$ 
matrices acting as linear operators on the complex Hilbert space 
$\bC^n$ with the standard inner product, $0_n$ resp.\ $\idty_n$ denoting
the zero resp.\ the multiplicative identity. We consider the Euclidean 
space of self-adjoint matrices endowed with the Hilbert-Schmidt inner 
product $(a,b)\mapsto\tr(ab)$ for $a,b\in{\rm Mat}(\bC,n)$ self-adjoint. 
Here $\tr$ denotes the standard trace. By $a\geq 0$ we mean that
$a\in{\rm Mat}(\bC,n)$ is positive semidefinite, i.e.\ self-adjoint and 
having non-negative eigenvalues. The {\it state space} of
${\rm Mat}(\bC,n)$ is the convex body
\begin{equation}
\label{defi:state_space}
\bS(n):=\{\rho\in{\rm Mat}(\bC,n)\mid\rho\geq 0\text{ and }\tr(\rho)=1\}.
\end{equation}
The Pauli 
$\sigma$-matrices 
$\sigma_1:=\left(\begin{smallmatrix}0&1\\1&0\end{smallmatrix}\right)$,
$\sigma_2:=\left(\begin{smallmatrix}0&-i\\i&0\end{smallmatrix}\right)$
and 
$\sigma_3:=\left(\begin{smallmatrix}1&0\\0&-1\end{smallmatrix}\right)$
together with $\idty_2$ are an orthogonal basis for the self-adjoint part 
of ${\rm Mat}(\bC,2)$. The {\it Bloch ball} is
\[\textstyle
\bS(2)=\{\frac{1}{2}(\idty_2+b_1\sigma_1+b_2\sigma_2+b_3\sigma_3)
\mid (b_1,b_2,b_3)\in \bB^3\}.
\]
For $m,n\in\bN$ the state space of the direct sum 
algebra $\cA:={\rm Mat}(\bC,m)\oplus{\rm Mat}(\bC,n)$ is the convex hull
of the individual state spaces
\[
\bS(\cA):=\bS(m+n)\cap\cA
=\conv(\bS(m)\oplus 0_n,0_m\oplus\bS(n)).
\] 
With $n$ direct summands we have e.g.\ the $n-1$ dimensional simplex
$\bS(\bC^n)$.
\par
An element $p\in\cA$ is an {\it orthogonal projection} if $p^2=p=p^*$. The 
set of orthogonal projections of $\cA$ are partially ordered by: $p\leq q$ 
if and only if $pq=p$ for $p,q$ orthogonal projections.
The {\it support projection} $s(\rho)$ of $\rho\in\bS(\cA)$ is the sum of 
the spectral projections of $\rho$ belonging to non-zero eigenvalues. The 
{\it maximal projection} $p_+(u)$ of a vector $u$ in the space 
$\cA_{\rm sa}$ of self-adjoint matrices is the spectral projection of $u$ 
for the largest eigenvalue of $u$. For non-zero $u\in\cA_{\rm sa}$ we have
the exposed faces (see Weis \cite{Weis}, Section 2.3)
\begin{equation}
\label{eq:s_faces}
\begin{array}{rcl}
F_{\perp}(\bS(\cA),u) & = &
 \{\rho\in\bS(\cA)\mid s(\rho)\leq p_+(u)\},\\
\ri(F_{\perp}(\bS(\cA),u)) & = & \{\rho\in\bS(\cA)\mid s(\rho)=p_+(u)\}
\end{array}
\end{equation}
and for $\rho\in\bS(\cA)$ we have the normal cones
\begin{equation}
\label{eq:s_cones}
\begin{array}{rcl}
\nc(\bS(\cA),\rho) & = & \{u\in\cA_{\rm sa}\mid s(\rho)\leq p_+(u)\},\\
\ri(\nc(\bS(\cA),\rho)) & = & \{u\in\cA_{\rm sa}\mid s(\rho)=p_+(u)\}.
\end{array}
\end{equation}
Much more general the facial structure of the state space of C*-algebra is 
treated by Alfsen and Schultz \cite{alfsen}.
It is immediate from (\ref{eq:s_faces}) and (\ref{eq:s_cones}) that every 
non-zero vector $u\in\cA_{\rm sa}$ is sharp normal for $\bS(\cA)$ and 
every element $\rho\in\bS(\cA)$ is sharp exposed in $\bS(\cA)$. We will 
extend this example in Example~\ref{ex:ip}.
\end{Exa}
\par
The definitions (\ref{eq:sharp_duality_1}) and (\ref{eq:sharp_duality_2})
depend {\it a priori} on the ambient space $\bE$ through the normal cone.
For sharp normal vectors we show independence in the following lemma. To 
keep notation clear we use orthogonal projections $\pi_V$ onto a vector 
space $V\subset\bE$ and not onto an affine space.
\begin{Lem}
\label{lem:gen_an}
Let $C\subset V$. Then every non-zero $v\in V^{\perp}$ is sharp normal for 
$C$ in the ambient space $\bE$. A vector $v\in\bE\setminus V^{\perp}$ is 
sharp normal for $C$ in the ambient space $\bE$ if and only if the vector 
$\pi_V(v)$ is sharp normal for $C$ in the ambient space $V$.
\end{Lem}
{\em Proof:\/}
For $v\in V^{\perp}\subset\lin(C)^{\perp}$ we have $F_{\perp}(C,v)=C$
(notice that $h(C,v)=0$ unless $C=\emptyset$). Then for every $x\in\ri(C)$
the normal cone $\nc(C,x)=\lin(C)^{\perp}$ is a vector space by 
Lemma~\ref{lin_complement}, so $v\in\ri(\nc(C,x))$ and $v$ is sharp normal 
for $C$.
\par
If $v\in\bE\setminus V^{\perp}$ then we have
$F_{\perp}(C,v)=F_{\perp}(C,\pi_V(v))$. For a point 
$x\in\ri(F_{\perp}(C,v))$ we distinguish between the normal cone
$\nc_\bE(C,x)$ in the ambient space $\bE$ and the normal cone
$\nc_V(C,x)\subset V$ in the ambient space $V$. These satisfy
$\nc_\bE(C,x)=\nc_V(C,x)+V^{\perp}$. By the sum formula (\ref{eq:risum})
for the relative interior we have
\[
\ri(\nc_\bE(C,x))=\ri(\nc_V(C,x))+V^{\perp}.
\]
Then we get $v\in\ri(\nc_\bE(C,x))$ if and only if 
$\pi_V(v)\in\ri(\nc_V(C,x))$, i.e.\ $v$ is sharp normal for $C$ in $\bE$ 
if and only if $\pi_V(v)$ is sharp normal for $C$ in $V$.
\hspace*{\fill}$\Box$\\
\par
Sharp normal vectors are preserved under projection.
\begin{Pro}
\label{pro:ac_projection}
If a non-zero vector $v\in V$ is sharp normal for $C$, then $v$ is sharp 
normal for $\pi_V(C)$. 
\end{Pro}
{\em Proof:\/}
We choose $x\in\ri(F_{\perp}(\pi_V(C),v))$ and we have to show that
$v\in\ri(\nc(\pi_V(C),x))$. By Lemma~\ref{lemma:lifted_faces} we have
\[
F_{\perp}(\pi_V(C),v)=\pi_V(F_{\perp}(C,v))
\]
so by (\ref{eq:cori}) we can choose a point $a\in\ri(F_{\perp}(C,v))$ such 
that $x=\pi_V(a)$. By assumption the vector $v$ is sharp normal for $C$ 
so $v\in\ri(\nc(C,a))$. By the formula for normal cones of a projected set
in Lemma~\ref{lemma:normal_cone_cylinder} we have
\[
\nc(\pi_V(C),x)=(\nc(C,a)\cap V)+V^{\perp}.
\]
Since $v\in\ri(\nc(C,a))$ we have $v\in\ri(\nc(C,a)\cap V)$ by the 
intersection formula (\ref{eq:riint}) for relative interiors. The sum 
formula (\ref{eq:risum}) for the relative interior shows
$v\in\ri(\nc(\pi_V(C),x))$, i.e.\ $v$ is sharp normal for $\pi_V(C)$ in 
$\bE$.
\hspace*{\fill}$\Box$\\
\par
We shortly discuss sharp exposed points and connect these to exposed faces.
The following lemma shows also that the definition 
(\ref{eq:sharp_duality_2}) of sharp exposed is independent of the ambient 
space because exposed faces are independent of the ambient space.
\begin{Lem}
\label{lem:sharp_exposed}
A non-empty face $F$ of $C$ is exposed if and only if there is a point in 
$\ri(F)$, which is sharp exposed in $C$. If there is a point in 
$\ri(F)$, which is sharp exposed in $C$, then all points in $\ri(F)$ are 
sharp exposed in $C$.
\end{Lem}
{\em Proof:\/}
Let $F$ be a non-empty exposed face of $C$. If $x\in\ri(F)$ then we have
$\nc(C,F)=\nc(C,x)$ by definition of the normal cone of $F$. We want to 
show that $x$ is sharp exposed in $C$. If $\nc(C,x)=\{0\}$ then there is 
nothing to prove. Otherwise by Lemma~\ref{lem:sup_exp} for all non-zero 
$u\in\ri(\nc(C,F))$ we have $F=F_{\perp}(C,u)$. In other words
for each $u\in\ri(\nc(C,x))\setminus\{0\}$ we have 
$x\in\ri(F_{\perp}(C,u))$, i.e.\ $x$ is sharp exposed in $C$.
\par
Conversely let $F\neq\emptyset$ be a face of $C$, not necessarily exposed. 
Since $C$ is exposed we can assume $F\neq C$, so $\nc(C,F)\neq\{0\}$ by
Lemma~\ref{lin_complement}. Let us choose a point $x\in\ri(F)$
and consider a non-zero vector $u\in\ri(\nc(C,F))=\ri(\nc(C,x))$. If we 
assume that $x$ is sharp exposed in $C$, then we have
$x\in\ri(F_{\perp}(C,u))$. Therefore $F=F_{\perp}(C,u)$ is an exposed face 
by the decomposition (\ref{eq:stratum}).
\hspace*{\fill}$\Box$\\
\par
Exposed faces are preserved under intersection.
\begin{Lem}
\label{lemma:intersection_all_exposed}
Let $\bA\subset\bE$ be an affine subspace and let $x\in C\cap\bA$. If  
$F(C,x)$ is an exposed face of $C$, then $F(C\cap\bA,x)$ is an exposed 
face of $C\cap\bA$.
\end{Lem}
{\em Proof:\/}
If $x\in\ri(C)$ then $x\in\ri(C\cap\bA)$ by the intersection formula 
(\ref{eq:riint}) for relative interiors. So $F(C\cap\bA,x)=C\cap\bA$ is 
exposed. Otherwise there is a non-zero $u\in\bE$ such that 
$x\in\ri(F_{\perp}(C,u))$. As $x\in\bA$ we have
$h(C,u)=\langle u,x\rangle=h(C\cap\bA,u)$, so we obtain
$F_{\perp}(C,u)\cap\bA=F_{\perp}(C\cap\bA,u)$. By the intersection formula 
(\ref{eq:riint}) for relative interiors this gives
$x\in\ri(F_{\perp}(C\cap\bA,u))$ and completes the proof.
\hspace*{\fill}$\Box$\\
%
%
%
%
%
%
%
\section{Touching cones}
\label{sec:tc}
\par
Let $C$ be a convex subset of a finite-dimensional real Euclidean vector 
space $(\bE,\langle\cdot,\cdot\rangle)$. We connect sharp normal vectors 
for $C$ to Schneider's \cite{schneider} concept of touching cone. Touching 
cones form a complete lattice with infimum the intersection. They include 
all normal cones, which are preserved under projection. Touching cones can 
detect the exposed faces which are intersections of coatoms.
\begin{Def}
\label{def:touching}
If $v\in\bE$ is a non-zero vector and if the exposed face 
$F_{\perp}(C,v)$ is non-empty, then the {\it touching cone} of $C$ for $u$ 
is defined by $\tc(C,u):=F(\nc(C,F_{\perp}(C,u)),u)$. This is the face of 
the normal cone $\nc(C,F_{\perp}(C,u))$, which has  $u$ in the relative 
interior. The normal cones $\lin(C)^\perp$ and $\bE$ are touching cones by 
definition, called {\it improper}. All other touching cones are
{\it proper}. The set of touching cones of $C$, called {\it touching cone 
lattice} is denoted by $\cT(C)$.
\end{Def}
\par
Perhaps the analogy with the face-function (as studied by
Klee and Martin \cite{klee} and others) should be pointed out here. The 
face-function associates with each $x\in C$ the smallest face $F(C,x)$ of 
$C$ containing $x$. Analogously (or dually) Definition~\ref{def:touching} 
associates with each vector $u\neq 0$ the smallest touching cone of $C$ 
containing it.
\begin{Lem}
\label{lemma:touching_cones_prop}
If $T$ is a touching cone of $C$, then $T=(T\cap\lin(C))+\lin(C)^\perp$.
Every normal cone of $C$ is a touching cone of $C$. If $T$ is a touching 
cone of $C$ but $T\neq\bE$, then
\[\begin{array}{ll}
\textup{(a)} &
\text{if } u\in\ri(T)\setminus\{0\}, \text{ then }
F_{\perp}(C,u)=\bigcap_{v\in T\setminus\{0\}}F_{\perp}(C,v)
\text{ is non-empty,}\\
\textup{(b)} &
\text{if } u\in\ri(T)\setminus\{0\}, \text{ then }
T=T(C,u),\\
\textup{(c)} &
\text{if } 0\in\ri(T), \text{ then } T=\lin(C)^{\perp}.
\end{array}\]
\end{Lem}
{\em Proof:\/}
The first assertion is clear for $T=\lin(C)^\perp$ or $T=\bE$. The normal 
cone $N$ of $x\in C$ is a direct sum of $N\cap\lin(C)$ and of
$\lin(C)^\perp$ by Lemma~\ref{lin_complement}, so this holds also for all 
its faces including $T$.
\par
Let us prove that every proper normal cone $N$ of $C$ belongs to $\cT(C)$. 
By the antitone lattice isomorphism $\cF_\perp(C)\to\cN(C)$ in 
Prop.~\ref{lattice_antitone_iso} there is a proper exposed face $F$, such 
that $N=\nc(C,F)$. By Lemma~\ref{lem:sup_exp} there exists 
$u\in\ri(\nc(C,F))\setminus\{0\}$ such that $F=F_{\perp}(C,u)$. Now 
$u\in\ri(\nc(C,F))=\ri(\nc(C,F_{\perp}(C,u)))$ gives $\tc(C,u)=\nc(C,F)$ 
by definition of a touching cone.
\par
(a)--(c) are trivial if $T=\{0\}$. Otherwise the touching cone $T$ arises 
from a non-zero vector $w\in\bE$ as $T=T(C,w)$ such that 
$F_{\perp}(C,w)\neq\emptyset$ (also in the case 
$T=\lin(C)^\perp\neq\{0\}$). 
\par
To show (a) we notice $T\subset\nc(C,F_{\perp}(C,w))$, so the intersection 
$\bigcap_{v\in T\setminus\{0\}}F_{\perp}(C,v)$ is non-empty by 
Lemma~\ref{lem:sup_exp}. For any $u\in\ri(T)\setminus\{0\}$ this 
intersection equals $F_{\perp}(C,u)$ by Cor.~\ref{cor:ex_inf}.
\par
To prove (b) we recall $w\in\ri(T)$ by definition of a touching cone. If a 
non-zero $u\in\ri(T)$ is chosen then by (a) we have 
$F_{\perp}(C,u)=F_{\perp}(C,w)$ and the two vectors $u,w$ belong to the 
relative interior $\ri(T)$ of the same face $T$ of $\nc(C,F_\perp(C,u))$,
so $T(C,u)=T(C,w)=T$ by the partition (\ref{eq:stratum}) of a convex set 
into relative interiors of faces.
\par
For (c) we recall that a convex cone with zero in the relative interior is 
a linear space. Since $w\in\ri(T)$ the opposite vector $-w$ belongs also 
to $\ri(T)$ and from (a) follows $F_{\perp}(C,w)=F_{\perp}(C,-w)$ so 
$F_{\perp}(C,w)=C$. The normal cone of $C$ is $\nc(C,C)=\lin(C)^{\perp}$ 
by Lemma~\ref{lin_complement} hence $T=T(C,w)=\lin(C)^{\perp}$.
\hspace*{\fill}$\Box$\\
\begin{Rem}
\label{rem:tc_cover}
If $K$ is a convex body and $u\in\bE$ a non-zero vector, then 
$F_\perp(K,u)$ is a non-empty exposed face and the touching cone 
$T:=T(C,u)$ with $u\in\ri(T)$ is defined. So $\bE\setminus\{0\}$ is 
covered by the relative interiors of touching cones $\neq\bE$.
Lemma~\ref{lemma:touching_cones_prop} (b) and (c) make sure that this 
cover is disjoint. We notice that this partition follows also from
Thm.~\ref{thm:lattice_isomo}
\end{Rem}
\par
Next we show beyond $\cT(C)\supset\cN(C)$ that the touching cone lattice 
consists of all non-empty faces of normal cones.  The infimum in $\cT(C)$ 
is the intersection and $\cT(C)$ is a complete lattice
\begin{equation}
(\cT(C),\subset,\cap,\vee,\lin(C)^\perp,\bE).
\end{equation}
\begin{Thm}
\label{thm:rtc_cover}
The touching cones of $C$ are exactly the non-empty faces of the normal 
cones of $C$, i.e.\ $\cT(C)=\{T\mid T\neq\emptyset \text{ is a face of } N,
N\in\cN(C)\}$. The touching cone lattice is a complete lattice ordered by 
inclusion. If $\{T_\alpha\}_{\alpha\in I}\subset\cT(C)$ is a non-empty 
family of touching cones, then
$\bigwedge_{\alpha\in I}T_\alpha=\bigcap_{\alpha\in I}T_\alpha$ and this 
intersection is a face of $T_{\widetilde{\alpha}}$ for every 
$\widetilde{\alpha}\in I$ with $T_{\widetilde{\alpha}}\neq\bE$.
\end{Thm}
{\em Proof:\/}
By definition, every touching cone is a non-empty face of a normal cone.
For the converse we need not treat the improper cones $\lin(C)^\perp$ and 
$\bE$, they have only one non-empty face, which is already included to the 
touching cones. Let $N$ be a proper normal cone of $C$. By the partition 
(\ref{eq:stratum}) of $N$ into relative interiors of its faces, it is 
sufficient to show for any non-zero vector $v\in N$ that $T(C,v)=F(N,v)$, 
i.e.\ the touching cone of $v$ is the face of $N$ with $v$ in the relative 
interior. 
\par
There is a proper exposed face $F$ with normal cone $\nc(C,F)=N$ by 
Prop.~\ref{lattice_antitone_iso}. Since $v\in\nc(C,F)$ we have
$F\subset F_{\perp}(C,v)$ as proved in Lemma~\ref{lem:face_cone_dual_u}.
By the antitone assignment of normal cones we get 
$\nc(C,F_\perp(C,v))\subset N$ and this statement includes by 
Prop.~\ref{normal_cone_intersection} that $\nc(C,F_\perp(C,v))$ is a face 
of $N$. By definition of a touching cone, $T(C,v)$ is a face of the normal 
cone $\nc(C,F_\perp(C,v))$, so it is a face of $N$.
As $v$ belongs to the relative interior of $T(C,v)$, we conclude that
$T(C,v)=F(N,v)$.
\par
In order to prove that $\cT(C)$ is a complete lattice with intersection as 
infimum, we can show by Remark~\ref{rem:closure_complete} for a non-empty 
family $\{T_\alpha\}_{\alpha\in I}$ that the intersection
$\bigcap_{\alpha\in I}T_\alpha$ is a touching cone of $C$. Since 
$\lin(C)^\perp$ is the smallest element of $\cT(C)$ by 
Lemma~\ref{lemma:touching_cones_prop} and since $\bE$ is the greatest 
element of $\cT(C)$ we assume that all $T_\alpha$ are proper touching 
cones. Then for every $\alpha\in I$ there is a non-zero $u_\alpha\in\bE$ 
such that $T_\alpha=T(C,u_\alpha)$. We put 
$N_\alpha:=\nc(C,F_\perp(C,u_\alpha))$ so $T(C,u_\alpha)$ is a face of
$N_\alpha$. The normal cone
$N:=\bigcap_{\widetilde{\alpha}\in I}N_{\widetilde{\alpha}}$ is a face 
of $N_\alpha$ by Prop.~\ref{normal_cone_intersection}, so the intersection
$N\cap T(C,u_\alpha)$ is a face of $N_\alpha$ and also of $N$. But then 
\[\textstyle
\bigcap_{\widetilde{\alpha}\in I}T_{\widetilde{\alpha}}
=\bigcap_{\widetilde{\alpha}\in I}(N\cap T(C,u_{\widetilde{\alpha}}))
\]
is a face of $N$, which is a touching cone by the first part of this 
theorem. Since the normal cone $N$ is a face of $N_\alpha$ the 
intersection $\bigcap_{\widetilde{\alpha}\in I}T_{\widetilde{\alpha}}$ is 
a face of $N_\alpha$.
\hspace*{\fill}$\Box$\\
\par
We prove an independence of touching cones.
\begin{Cor}
\label{cor:ind_t}
The lattice orderings of $\cN(C)$ and $\cT(C)$ and the embedding
$\cN(C)\to\cT(C)$ are independent of the ambient space $\bE$.
\end{Cor}
{\em Proof:\/}
A normal cone $N\in\cN(C)$ has the direct sum form 
$N=(N\cap\lin(C))+\lin(C)^\perp$ by Lemma~\ref{lin_complement}. Thus the 
normal cone lattice $\cN(C)$ can be reconstructed from
\[
\widetilde{\cN}:=\{N\cap\lin(C)\mid N\in\cN(C)\}
\]
by adding the direct summand $\lin(C)^\perp$. This defines a lattice 
isomorphism $\widetilde{\cN}\to\cN(C)$ and the lattice $\widetilde{\cN}$ 
is independent of the ambient space $\bE$ because $\widetilde{\cN}$ is the 
normal cone lattice of $C$ in the ambient space $\lin(C)$.
By Thm.~\ref{thm:rtc_cover} the touching cone lattice $\cT(C)$ consists of 
all non-empty faces $T$ of $\cN(C)$, so $T=(T\cap\lin(C))+\lin(C)^\perp$ 
holds. The same argument as above shows independence of the lattice 
$\cT(C)$ from the ambient space $\bE$. The question which touching cones 
are normal cones is solved by the embedding $\cN(C)\to\cT(C)$, which is 
also induced from the ambient space $\lin(C)$.
\hspace*{\fill}$\Box$\\

\par
Sharp normal vectors characterize the normal cones among all touching 
cones.
\begin{Pro}
\label{prop:sharp_normal}
A proper touching cone $T$ of $C$ is a normal cone of $C$ if and only if 
there is a vector in $\ri(T)\setminus\{0\}$, which is sharp normal for 
$C$. If there is a vector in $\ri(T)\setminus\{0\}$, which is 
sharp normal for $C$, then all vectors in $\ri(T)\setminus\{0\}$ are sharp 
normal for $C$.
\end{Pro}
{\em Proof:\/}
Let $K$ be a proper touching cone of $C$ and let us assume that 
$u\in\ri(K)\setminus\{0\}$ is sharp normal for $C$. Then there exists
$x\in\ri(F_{\perp}(C,u))$ and we have $u\in\ri(\nc(C,x))$. By definition 
of the normal cone of a face we have $\nc(C,x)=\nc(C,F_{\perp}(C,u))$ 
hence $u\in\ri(\nc(C,F_{\perp}(C,u)))$ and this gives us
$T(C,u)=\nc(C,F_{\perp}(C,u))$. Since $u\in\ri(K)$ we have $K=T(C,u)$ by
Lemma~\ref{lemma:touching_cones_prop} (b). Hence $K$ is the 
normal cone of the non-empty face $F_{\perp}(C,u)$.
\par
Conversely let us assume that the touching cone $K$ is the normal 
cone of a non-empty face of $C$. Then by 
Prop.~\ref{lattice_antitone_iso} we have $K=\nc(C,F)$ for some
non-empty exposed face $F$ of $C$. Now Lemma~\ref{lem:sup_exp} shows 
for any non-zero $u\in\ri(K)$ that $F=F_{\perp}(C,u)$ holds. Then for any
$x\in\ri(F_{\perp}(C,u))$ we have
\[
\nc(C,x)=\nc(C,F)=K
\]
and this shows that $u\in\ri(\nc(C,x))$. We have proved that $u$ is sharp 
normal for $C$. If $K$ is proper, then existence of a non-zero vector $u$ 
in $\ri(K)$ is assured.
\hspace*{\fill}$\Box$\\
\par
Projection properties of sharp normal vectors apply to touching cones. We 
denote $\pi_V$ the orthogonal projection onto a vector space $V\subset\bE$.
\begin{Cor}
\label{cor:trans_sn_tc}
Let $v\in V\setminus\{0\}$. If the touching cone $T(\pi_V(C),v)$ exists 
and is not a normal cone, then $T(C,v)$ exists and is not a normal cone.
In particular, if $\cT(C)=\cN(C)$ then $\cT(\pi_V(C))=\cN(\pi_V(C))$.
\end{Cor}
{\em Proof:\/}
If $T(\pi_V(C),v)$ exists, then $F_\perp(\pi_V(C),v)\neq\emptyset$ and
by Lemma~\ref{lemma:lifted_faces} we have
$F_{\perp}(\pi_V(C),v)=\pi_V(F_{\perp}(C,v))$. So 
$F_{\perp}(C,v)\neq\emptyset$ and $T(C,v)$ exists. If in addition $T(C,v)$ 
is a normal cone of $C$, then $v$ is sharp normal for $C$ by 
Prop.~\ref{prop:sharp_normal} as $v\in\ri(T(C,v))$. Then by
Prop.~\ref{pro:ac_projection} $v$ is sharp normal for $\pi_V(C)$ and this 
implies that $T(\pi_V(C),v)$ is a normal cone of $\pi_V(C)$.
\hspace*{\fill}$\Box$\\
\begin{Exa}
\label{ex:ip}
We return to Example~\ref{ex:state_spaces} and denote by $K:=\bS(\cA)$ the 
state space of the algebra $\cA:={\rm Mat}(\bC,2)\oplus\bC$. We have seen 
that every non-zero $u\in\cA_{\rm sa}$ is sharp normal for $K$ and every 
$\rho\in K$ is sharp exposed in $K$. This implies $\cT(K)=\cN(K)$ and 
$\cF(K)=\cF_\perp(K)$ by Prop.~\ref{prop:sharp_normal} and 
Lemma~\ref{lem:sharp_exposed}. Now we consider a family of 
two-dimensional projections and intersections of $K$ produced by a 
three-dimensional affine space of self-adjoint matrices without
$\sigma_3$-contribution in the first summand
\[
\widetilde{\bA}:=
\{a\in{\rm Mat}(\bC,2)\oplus\bC\mid a^*=a,
\tr[a(\sigma_3\oplus 0)]=0\text{ and }\tr(a)=1\}.
\]
If $\pi_{\widetilde{\bA}}$ denotes orthogonal projection to 
$\widetilde{\bA}$, then (see \cite{knauf_weis} Section 2)
\[
C:=\pi_{\widetilde{\bA}}(K)=K\cap\widetilde{\bA}=\conv\,
\left[\{\rho\in\bS(2)\mid\tr(\rho\sigma_3)=0\}\oplus 0,\,0_2\oplus 1\right]
\]
is the three-dimensional cone depicted in
Figure~\ref{fig:Kegel_Projektion}, left. 
By Coro.~\ref{cor:trans_sn_tc} we have $\cT(C)=\cN(C)$ because $C$ is the
projection of $K$ to $\widetilde{\bA}$. By 
Lemma~\ref{lemma:intersection_all_exposed} we have $\cF(C)=\cF_\perp(C)$ 
because $C$ is the intersection of $K$ with $\widetilde{\bA}$.
\par
Let $\bA\subset\widetilde{\bA}$ be the two-dimensional affine subspace 
containing $\frac{\idty_3}{3}$ and having the angle $\varphi$ with the
direction $-\idty_2\oplus 2$. Two example are shown in 
Figure~\ref{fig:Kegel_Projektion}, right. The projection shapes 
$\pi_{\bA}(C)$ have every touching cone a normal cone. So, according to 
Thm.~\ref{thm:local} and Remark~\ref{ref:exposed_2D} a face of 
$\pi_{\bA}(C)$ is non-exposed if and only if it is the endpoint of a 
unique one-dimensional face. The examples with 
$\varphi=12^\circ$ and $\varphi\approx39^\circ$ have two non-exposed
faces: the tangent points of boundary segments to the elliptic boundary 
arcs. The intersections $C\cap\bA$ have all faces exposed. In the depicted 
examples exist touching cones, which are not normal cones. It is 
instructive to 
realize that projection and intersection for the same affine space $\bA$
are polars of each other up to the sign (see e.g.\ Weis \cite{Weis} 
Section 2.4).
\end{Exa}
\begin{figure}[t!]
\centerline{%
\begin{picture}(14,3.5)
\ifthenelse{\equal{\ftype}{eps}}{%
\put(0,0){\includegraphics[height=3.5cm,bb=0 0 400 230]{Kegel_Aspekt.eps}}
\put(6.2,0){\includegraphics[height=3.5cm,bb=0 50 400 300]{%
Proion12.eps}}
\put(10,0){\includegraphics[height=3.5cm,bb=0 0 400 250]{%
Proion39.eps}}}{%
\put(0,0){\includegraphics[height=3.5cm,bb=0 0 400 285]{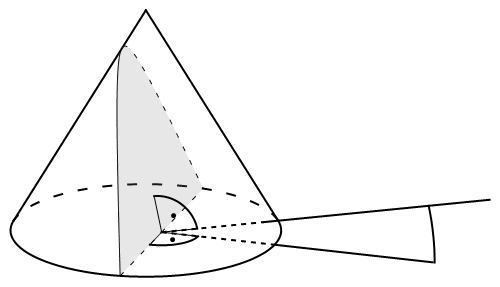}}
\put(6.7,0){\includegraphics[height=3.5cm,bb=30 50 400 360]{%
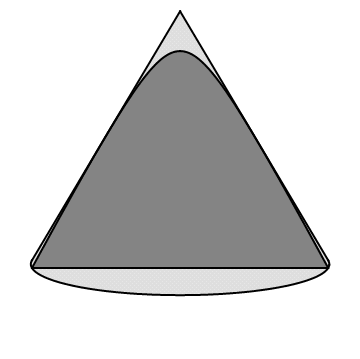}}
\put(12,0){\includegraphics[height=3.5cm,bb=170 -5 400 305]{%
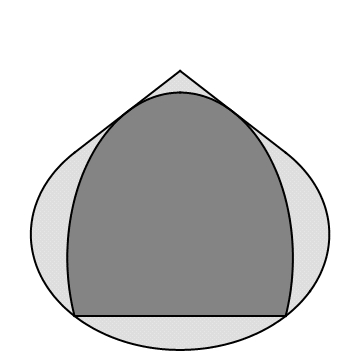}}}
\put(6.5,3){$\varphi=12^\circ$}
\put(10,3){$\varphi\approx 39^\circ$}
\put(1.6,2){$\bA$}
\put(4.7,0.6){$\varphi$}
\end{picture}}
\caption{%
\label{fig:Kegel_Projektion}The cone of revolution $C$ of an equilateral 
triangle (left). An affine plane $\bA$ through the center of gravity of $C$
is specified by the angle $\varphi<90^\circ$. The intersection $C\cap \bA$ 
is hyperbolic for $\varphi=12^\circ$ (middle) and elliptic for 
$\varphi\approx39^\circ$ (right), drawn dark. The union with the 
bright region surrounding it is the projection $\pi_{\bA}(C)$.}
\end{figure}
\par
An easy corollary of Minkowski's and Carath\'eodory's theorem  
characterizes normal cones and exposed faces in terms of touching cones.
\begin{Thm}
\label{thm:local_n}
Let $N$ be a proper normal cone of $C$ such that every touching cone
included in $N$ is a normal cone. Then $N$ can be written as a supremum of 
atoms of $\cN(C)$. A number of $\dim(N)-\dim(\lin(C)^\perp)$ atoms suffice 
in the supremum.
\end{Thm}
{\em Proof:\/}
By Coro.~\ref{cor:ind_t} we assume that $C$ has non-empty interior
${\rm int}(C)\neq\emptyset$, so $\{0\}$ is the smallest element 
in $\cN(C)$. Let $N$ be a proper normal cone of $C$. Then $N$ does not 
contain a line, for otherwise by (iv) in (\ref{eq:nc_for_begn}) we had 
${\rm int}(C)=\emptyset$. Therefore there is an affine hyperplane 
$H\subset\bE$ such that $K:=N\cap H$ is a convex body and $N$ is the 
positive hull $N=\pos(K)$. Let $u\in\ri(K)$. By Minkowski's theorem we 
write $u$ as a convex combination $u=\sum_{i=1}^d\lambda_iu_i$ for 
(non-zero) extreme points $u_i$ of $K$. By Carath\'eodory's theorem we 
choose $d=\dim(K)+1=\dim(N)$.
\par
We show that $N$ is a supremum of the $d$ normal cones 
$r_i:=\{\lambda u_i\mid\lambda\geq 0\}$, $i=1,\ldots,d$. By 
Lemma~\ref{lem:cone_faces} $u$ belongs to $\ri(N)$ and the ray $r_i$ is a 
face of $N$ so $r_i$ is a touching cone by Thm.~\ref{thm:rtc_cover}. By 
assumption the touching cone $r_i$ is a normal cone so it is an atom in 
$\cN(C)$. If the supremum $\widetilde{N}:=r_1\vee\cdots\vee r_d$ is 
strictly included into $N$, then $\widetilde{N}$ must be a proper face of 
$N$ by Prop.~\ref{normal_cone_intersection}, so 
$\widetilde{N}\subset\rb(N)$ by the partition (\ref{eq:stratum}) of $N$ 
into relative interiors of its faces. This contradicts $u\in\ri(N)$.
\hspace*{\fill}$\Box$\\
\par
The isomorphism 
$\nc(C):\cF_{\perp}(C)\to\cN(C)$ in Prop.~\ref{lattice_antitone_iso} gives 
an equivalent form of this theorem (which is trivial if $C$ is a single 
point). 
\begin{Thm}
\label{thm:local}
Let $F$ be a proper exposed face of $C$ such that every touching cone
included in the normal cone $\nc(C,F)$ is a normal cone. Then $F$ can be 
written as an intersection of coatoms of $\cF_\perp(C)$. A number of 
$\dim(\nc(C,F))-\dim(\lin(C)^\perp)$ coatoms suffice in the intersection.
\end{Thm}
\par
One may check Thm.~\ref{thm:local_n} and 
Thm.~\ref{thm:local} on
Figure~\ref{fig:sharp_non-closed},~\ref{fig:lattices}
and~\ref{fig:quarter_disk}. The theorems have no converse by example in
Figure~\ref{fig:lens}. The bound on coatoms is saturated by a corner of a 
cube, it is not saturated for the apex of the cone in 
Figure~\ref{fig:Kegel_Projektion}, left.
%
%
%
%
%
%
\section{Polar convex bodies}
\label{sec:polarity}
\par
This section is restricted to a convex body $K\subset\bE$ in a 
finite-dimensional real Euclidean vector space 
$(\bE,\langle\cdot,\cdot\rangle)$. Unless specified other we assume that 
$K$ has non-empty interior ${\rm int}(K)\neq\emptyset$ containing the 
origin $0\in{\rm int}(K)$ and second we assume that $K$ has at least two 
points. Conjugate faces induce an isotone lattice isomorphism between the 
faces of the polar convex body $K^\circ$ and the touching cones of $K$. 
This implies an equivalent theorem to Thm.~\ref{thm:local_n}, which can be 
proved directly using only Minkowski's and Carath\'eodory's theorem. The 
antitone lattice isomorphism $\cF_\perp(K)\to\cN(K)$ (see
Prop.~\ref{lattice_antitone_iso}) gives a fourth equivalent form of 
Thm.~\ref{thm:local_n}.
\begin{Def}
The {\it polar body} of $K$ is $K^{\circ}:=\{x\in\bE\mid\langle x,y\rangle
\leq 1\text{ for all }y\in K\}$. If $F$ is a subset of $K$, then the
{\it conjugate face} of $F$ is $\widehat{F}:=\{x\in K^\circ\mid
\langle x,y\rangle=1\text{ for all }y\in F\}$.
\end{Def}
\par
The polar body $K^\circ$ is a convex body with $0\in{\rm int}(K^\circ)$ 
and such that $K^{\circ\circ}=K$, see Schneider \cite{schneider}, Section 
1.6. An example of a convex body with its polar body is depicted in 
Figure~\ref{fig:exposed_faces_normal_cones}, right. We recall that 
$\emptyset$ and $K$ are exposed faces of $K$ so as to make $\cF_\perp(K)$ 
a lattice (this deviates from definitions by Rockafellar or Schneider 
\cite{rock,schneider}). By Schneider, Thm.~2.1.4, 
a subset $F\subset K$ is included in a proper exposed face of $K$ if and 
only if the conjugate face $\widehat{F}$ is a proper exposed face of 
$K^\circ$. Further, if these conditions hold, then 
$(\widehat{F})\,\raisebox{0.3ex}{$\widehat{}$}=\sex(F)$ is the smallest 
exposed face of $K$ containing $F$. Obviously 
$\widehat{\emptyset}=K^\circ$ and $\widehat{K}=\emptyset$. So 
$\cF_\perp(K)\to\cF_\perp(K)$, 
$F\mapsto(\widehat{F})\,\raisebox{0.3ex}{$\widehat{}$}\,$\, is the 
identity and we get an antitone lattice isomorphism:
\begin{equation}
\label{eq:conjugate}
\cF_\perp(K)\to\cF_\perp(K^\circ),\quad F\mapsto\widehat{F}.
\end{equation}
An example is shown in Figure~\ref{fig:lattices}. The following 
remark may help our intuition.
\begin{Rem}
The {\it polar} of an affine space $\bA$ in $\bE$ with respect to the unit 
sphere $\{x\in\bE\mid\langle x,x\rangle=1\}$ is the affine space
\[
\bA^{\rm polar}:=\{x\in\bE\mid\langle x,y\rangle=1
\text{ for all }y\in\bA\}.
\]
The polar is well-known in projective geometry (see e.g.\ Coxeter or 
Fischer \cite{coxeter,fischer}), it defines an antitone lattice 
isomorphism on the set of affine subspaces of $\bA\subset\bE$ with 
$0\not\in\bA$ with $\bE$ joined. The polar is an involution, i.e.\
$\bA^{\rm polar}{}^{\rm polar}=\bA$ such that
$\dim(\bA)+\dim(\bA^{\rm polar})=\dim(\bE)-1$. E.g.\
$\emptyset^{\rm polar}=\bE$ and $\bE^{\rm polar}=\emptyset$.
In fact it restricts a {\it correlation} of a projective space.
\par
The conjugate face of an arbitrary subset $F\subset K$ is 
$\widehat{F}=\aff(F)^{\rm polar}\cap K^\circ$. It is possible, e.g.\ for a 
disk, that $\aff(\widehat{F})\subsetneq\aff(F)^{\rm polar}$. Equality 
holds for all polytopes $K$ and their faces $F$, see Gr\"unbaum 
\cite{gruenbaum}, Section 3.4.
\end{Rem}
\par
The next observation is that the normal cone of every non-empty exposed 
face $F$ of $K$ is the positive hull of the conjugate face 
$\nc(K,F)=\pos(\widehat{F})$ (we have $\pos(\emptyset)=\{0\}$). This 
statement is proved in a more general form by Schneider \cite{schneider}, 
Lemma 2.2.3. We include the empty face with $\widehat{\emptyset}=K^\circ$ 
and with normal cone $\nc(K,\emptyset)=\pos(K^\circ)=\bE$. Combining this 
with the two antitone lattice isomorphisms 
$\cF_\perp(K)\to\cF_\perp(K^\circ)$ in (\ref{eq:conjugate}) and  
$\cF_\perp(K)\to\cN(K)$ in Prop.~\ref{lattice_antitone_iso} we get an 
isotone lattice isomorphism
\begin{equation}
\label{eq:cb_iso}
\cF_\perp(K^\circ)\to\cN(K),\quad F\mapsto\pos(F)
\end{equation}
from the commuting diagram
\[\left.\xymatrix{%
\cF_\perp(K) \ar[r] \ar[d] & \cF_\perp(K^\circ) \ar[dl]^{\pos}\\
\cN(K)
}\right..
\]
Every proper exposed face $F$ of $K^\circ$ has a supporting hyperplane $H$ 
of $K^\circ$ with $F=K^\circ\cap H$. We get $F=\pos(F)\cap H$ and since 
$\emptyset\neq{\rm int}(K^\circ)$ we have also
$F=\pos(F)\cap \rm(K^\circ)$. So the inverse to (\ref{eq:cb_iso}) is
\[
\cN(K)\to\cF_\perp(K^\circ),\quad
\left\{\begin{array}{lcl}
N\mapsto\rb(K^\circ)\cap N & \text{if} & N\neq\bE\\
\bE\mapsto K^\circ
\end{array}\right..
\]
By examples in Figure~\ref{fig:lattices} the antitone isomorphism 
$\cF_{\perp}(K)\to\cN(K)$ does not extend to $\cF(K)\to\cT(K)$ but we 
prove extension of $\pos:\cF_\perp(K^\circ)\longrightarrow\cN(K)$.
\begin{Thm}
\label{thm:lattice_isomo}
Let $K$ be a convex body containing at least two points and with
$0\in{\rm int}(K)$. If $K^\circ$ denotes the polar body, then the positive 
hull operator $\pos$ defines an isotone lattice isomorphism
$\cF(K^\circ)\to\cT(K)$.
\end{Thm}
{\em Proof:\/}
We consider a proper exposed face $F$ of $K^\circ$. For $N:=\pos(F)$ we 
have a bijection $\pos: \cF(F)\to\cF(N)\setminus\{\emptyset\}$ by 
Lemma~\ref{lem:cone_faces}, which may be written in the form
\begin{equation}
\label{eq:local_iso}
\begin{array}{lcl}
\pos(\widetilde{F})\cap F=\widetilde{F}
& \text{for all} & \widetilde{F}\in\cF(F),\\
\pos(G\cap F)=G
& \text{for all} & G\in\cF(N)\setminus\{\emptyset\}.
\end{array}\end{equation}
By (\ref{eq:cb_iso}) and the paragraph following it, we have 
$F=\rb(K^\circ)\cap N$ so we replace $F$ by $\rb(K^\circ)$ in 
(\ref{eq:local_iso}) except $\cF(F)$, which we leave unchanged. This gives 
us the bijection
\[
\pos:
\left\{\parbox{3.3cm}{faces of proper\\exposed faces of $K^\circ$}\right\}
\to
\left\{\parbox{4.4cm}{non-empty faces of proper\\normal cones of 
$K$}\right\}.
\]
The domain is clearly $\cF(K^\circ)\setminus\{K^\circ)$ and the target is 
$\cT(K)\setminus\{\bE)$ by Thm.~\ref{thm:rtc_cover}. Since $K$ has more 
than two points we have $\bE\neq\{0\}$ so $\pos(K^\circ)=\bE$ extends this 
map to an isotone lattice isomorphism $\cF(K^\circ)\to\cT(K)$.
\hspace*{\fill}$\Box$\\
\par
Theorem~\ref{thm:lattice_isomo} partitions $\bE\setminus\{0\}$ into 
relative interiors of touching cones, see Rem.~\ref{rem:tc_cover}.
We translate Thm.~\ref{thm:local_n} by interchanging exposed faces with 
normal cones and touching cones with faces, using (\ref{eq:cb_iso}) and 
Thm.~\ref{thm:lattice_isomo}. Through affine embeddings we can drop the 
condition $0\in\ri(K)$ in the sequel, the condition that $K$ has at least 
two points is not needed.
\begin{Thm}
\label{thm:minkowski}
Let $K$ be a convex body and let $F$ be a proper exposed face of $K$ such 
that every face included in $F$ belongs to $\cF_\perp(K)$. Then $F$ can be 
written as a supremum of at most $\dim(F)+1$ atoms of $\cF_\perp(K)$.
\end{Thm}
\par
Thm.~\ref{thm:minkowski} follows directly from Minkowski's and
Carath\'eodory's theorem. It is wrong if $K$ is not closed (e.g.\ a closed 
triangle with an extreme point missing) or unbounded (e.g.\ the strip 
$\{(x,y)\in\bR^2\mid x,y\geq 0, y\leq 1\}$).
\par
The antitone lattice isomorphism $\nc(K):\cF_{\perp}(K)\to\cN(K)$ in
Prop.~\ref{lattice_antitone_iso} gives us an equivalent form of 
Thm.~\ref{thm:minkowski}. We denote by $F_\perp(K,N)$ the unique exposed 
face $F$ of $K$ with $\nc(K,F)=N$ and we use intersection for the infimum 
in $\cN(K)$ by Prop.~\ref{normal_cone_intersection}. 
\begin{Thm}
Let $K$ be a convex body and let $N$ be a proper normal cone of $K$ such 
that every face included in $F_\perp(K,N)$ belongs to $\cF_\perp(K)$. Then 
$N$ can be written as an intersection of at most $\dim(F_\perp(K,N))+1$
coatoms of $\cN(K)$.
\end{Thm}
\par
The bound on coatoms is saturated by the normal vector of a square face
of the cube.
\vspace{0.5cm}
%
%
%
%
\par\noindent
\textbf{Acknowledgment.}
Thanks to Andreas Knauf for discussions and encouragement in the early 
stage of this work, to Hermann Schulz-Baldes for suggestions to improve 
the text, to Rolf Schneider for drawing my attention to the 
Aleksandrov-Fenchel inequality and for correcting errors and to the 
anonymous referee's useful comments.
%
%
%
%
%
\bibliographystyle{amsalpha}

\end{document}